\documentclass[12pt]{article}
\usepackage{amsmath}
\usepackage{amssymb}
\newcommand{\Alex}{{\cal A}_C}
\newcommand{\fhat}{\widehat{f}}

\newcommand{\Rbar}{\overline{\R}}
\newcommand{\Ct}{{\cal B}_C}
\newcommand{\D}{{\cal D}}
\newcommand{\Dp}{{\cal D}'}

\newcommand{\ebv}{{\cal EBV}}
\newcommand{\bv}{{\cal BV}}
\newcommand{\intinf}{\int^\infty_{-\infty}}

\newcommand{\N}{{\mathbb N}}
\newcommand{\R}{{\mathbb R}}
\newcommand{\C}{{\mathbb C}}
\newcommand{\fn}{\!:\!}
\newcommand{\lsum}{\sum\limits}

\newcommand{\qed}{\mbox{$\quad\blacksquare$}}
\newtheorem{theorem}{Theorem}

\newtheorem{prop}[theorem]{Proposition}
\newtheorem{corollary}[theorem]{Corollary}

\newtheorem{defn}[theorem]{Definition}
\newtheorem{example}[theorem]{Example}
\begin{document}
\hspace{-2cm}
\raisebox{12ex}[1ex]{\fbox{{\footnotesize
Preprint
December 7, 2007.\quad To appear in {\it Real Analysis Exchange.}
}}}

\begin{center}
{\large\bf The distributional Denjoy  integral}
\vskip.25in
Erik Talvila\footnote{Supported by the
Natural Sciences and Engineering Research Council of Canada.
}\\ [2mm]
{\footnotesize
Department of Mathematics and Statistics \\
University College of the Fraser Valley\\
Abbotsford, BC Canada V2S 7M8\\
Erik.Talvila@ucfv.ca}
\end{center}

{\footnotesize
\noindent
{\bf Abstract.} 
Let $f$ be a distribution (generalised function) on the real line.  If
there is a continuous function $F$ with real limits at infinity such that
$F'=f$ (distributional derivative) then
the distributional integral of $f$ is defined as $\int_{-\infty}^\infty f =
F(\infty) - F(-\infty)$.  It is shown that this simple definition gives
an integral that includes the Lebesgue and Henstock--Kurzweil integrals.
The Alexiewicz norm leads to a Banach space of integrable distributions
that is isometrically isomorphic to the space of 
continuous functions on the extended
real line with uniform norm.  The
dual space is identified with the functions of bounded variation.
Basic properties of integrals are established using elementary properties
of distributions: integration by parts, H\"older inequality, 
change of variables, convergence theorems, Banach lattice structure, 
Hake theorem, Taylor theorem, second mean value theorem.  Applications are
made  to
the half plane Poisson integral and Laplace transform.  The paper 
includes a short history of
Denjoy's descriptive integral definitions.  Distributional integrals
in Euclidean spaces are discussed and a more general distributional
integral that also integrates Radon measures is proposed.
\\
{\bf 2000 subject classification:} 26A39, 46B42, 46E15, 46F05, 46G12\\
{\bf Key words:} distributional Denjoy integral; continuous primitive
integral; Henstock--Kurzweil integral; 
Schwartz distributions; Alexiewicz norm; Banach lattice.
}\\

\section{Introduction}

We are fortunate to live in a richly diverse universe in which there are
many integrals and many interesting ways of defining these integrals.
Some of the major integrals are those of Riemann, Lebesgue, Denjoy and
Henstock--Kurzweil.  
In this paper we will present a theory of integration based on the 
descriptive Denjoy method. 
The definition is simple and elegant.  A distribution $f$ is integrable if there
is a continuous function $F$ whose distributional derivative equals
$f$.  Then $\int_a^bf=F(b)-F(a)$.    
This is a very powerful integral that includes all of those
mentioned above. To define it we only need the notion of distributional
derivative and the Riemann  integration of continuous functions.  
No measure theory is needed to define the integral 
and there are no partitions to construct.
We will see that under the Alexiewicz norm (see Section~\ref{arnauddenjoy}),
 the
space of integrable distributions forms a Banach
space (and Banach lattice) that is isometrically isomorphic
to the space of continuous functions on the extended real
line with uniform norm.
The dual space is identified with the space of functions of bounded
variation.
There are general versions of the Fundamental Theorem of
Calculus,
integration by parts and change
of variables formulas, a H\"older inequality, convergence
theorems,  Taylor's theorem, Hake's theorem and
the second mean value theorem.  We give applications to the half plane
Poisson integral and the Laplace transform.  Absolute integration is also discussed.  All of these results are easy
to prove using only elementary results in distributions 
(generalised functions). Besides distributions,  
we will assume some familiarity with Riemann--Stieltjes integrals,
functions of bounded variation, and basic notions of
functional analysis, such as Cauchy sequences in the Banach space
of continuous functions with uniform norm $\|\cdot\|_\infty$.
Most of the results we use in distributions are summarised in
Section~\ref{distributions}.
The reader should have a
nodding acquaintance with Lebesgue measure
and integration although it will be apparent that this approach to
integration de-emphasises measure and puts more emphasis on functional
analytic aspects.
Our setting will be integration on the
real line with respect to Lebesgue measure $\lambda$.
At the end of
the paper we sketch out generalisations to integration in $\R^n$ and
integration  with
respect to Radon measures.

\section{Integrating derivatives}\label{arnauddenjoy}

The Riemann and Lebesgue integrals are both absolute.  This means that if
function
$f$ is integrable, so is $|f|$.  An outcome of this is that we get a weaker
version of the Fundamental Theorem than we'd like.  For example, the
function $F(x)=x^2\cos(x^{-2})$
for $x\not=0$ and $F(0)=0$ is differentiable at each point in $\R$ but $F'$
is not continuous at $0$ since $F'(x)\sim 2x^{-1}\sin(x^{-2})$ as $x\to 0$.  
And, $\int_0^1F'$ does not exist in the Riemann
or Lebesgue sense since $|F'|$ is not integrable in a neighbourhood of $0$.
However,
$\int_0^1F'$ exists as a conditionally convergent improper Riemann integral
and hence as a Henstock--Kurzweil integral.  So, to be able to write
$\int_0^1F'=F(1)-F(0)$ we need to consider nonabsolute integrals.
The 
problem of integrating derivatives was solved 
by Denjoy in the early part of the 20th century.

Arnaud Denjoy (pronounced rather like ``dawn-djwah'') was a French mathematician
who was born in 1884 and lived for over 90 years.   He produced three different solutions
to the problem of integrating derivatives and is known for several other 
results in function  theory, Fourier series, quasi-analytic functions and
dynamical systems.  See \cite{mactutor} for a photo.

Denjoy's solution was to use a descriptive definition of the integral.
This defines the integral via its {\em primitive}.  This is a
continuous function whose derivative in some sense is equal to the integrand.
For example, $F\fn\R\to\R$ is {\em absolutely continuous} ($AC$) if
for every $\epsilon>0$ there is $\delta>0$ such that whenever
$\{(x_n,y_n)\}$ is a sequence of disjoint intervals with $\sum|x_n-y_n|<\delta$
we have $\sum|F(x_n)-F(y_n)|<\epsilon$.
This definition readily generalises to
arbitrary measure spaces.  We have the strict inclusions
$C^1\subsetneq AC\subsetneq C^0$.  If $F$ is $AC$ then it is continuous
on $\R$ and is differentiable almost everywhere. 
If $F\in AC$ then $\int_a^bF'=F(b)-F(a)$.
The descriptive definition of the Lebesgue integral is then:
$f\fn[a,b]\to\R$ is integrable if there is a function
$F\in AC$, called the primitive, such
that $F'=f$ almost everywhere.   In this case, $\int_a^bf=F(b)-F(a)$.
This is one half of the Fundamental
Theorem of Calculus.  The other half says that if $f\in L^1$ then
$F(x):=\int_a^xf$ defines an $AC$ function and $F'=f$
almost everywhere.  The function $F(x)=x^{2}\cos(x^{-2})$ at the
beginning of this section is not $AC$.

The corresponding function space for Denjoy  integrals is
$ACG*$ (generalised absolute continuity in the restricted sense).  
The precise definition need not concern us here.  If  you are 
interested, see \cite{gordon}.  
The important thing is 
that $C^1\subsetneq AC\subsetneq ACG*\subsetneq C^0$ and we have
a larger, more complicated space in which functions have derivatives
almost everywhere.  The Denjoy integral is then defined by saying
that $f$ is integrable if it has a primitive $F\in ACG*$ such that
$F'=f$ almost everywhere.  Then, $\int_a^bf=F(b)-F(a)$. 

Since $AC\subsetneq ACG*$, the Denjoy integral properly contains the Lebesgue
integral (with respect to Lebesgue measure on the real line). 
It turns out that if a continuous function is differentiable
everywhere then it is in $ACG*$.  The same is true if the function has a
derivative everywhere except in a countable set.  Hence, we can integrate
the function $F'$ given at the beginning of this section.  
The Denjoy integral is equivalent to the Henstock--Kurzweil integral,
which is defined using Riemann sums.  It is also equivalent to the
Perron integral, which is defined using major and minor functions \cite{gordon}.

The Denjoy integrable functions are made into a normed linear space
via the Alexiewicz norm \cite{alexiewicz}.  
This is defined by $\|f\|=\sup_{a\leq x\leq b}|\int_a^xf|$.
Unfortunately, this does not define a Banach space so we do not have
analogues of the many wonderful results in $L^p$ spaces.  
Real analysts delight in working with spaces such as $ACG*$ (see
any volume of the journal {\it Real Analysis Exchange}). However,
the attraction of 
such spaces has been less compelling for other mathematicians.
One problem is that there is no canonical generalisation to $\R^n$.
A considerable amount of research was carried out in Denjoy integration
until the end of the 1930's but these
deficiencies caused this integral to be virtually abandoned by 1940.  
However, we get a much simpler and
yet more powerful integral by using the distributional Denjoy integral.
For this, we will need to briefly introduce some results in distributions.

\section{Schwartz distributions}\label{distributions}

The theory of distributions, or generalised functions, extends the
notion of function so that we no longer have pointwise values but
all distributions have derivatives of all orders.  Most of the final
theory that emerged in the 1940's was due to Laurent Schwartz but 
of course he did not work in vacuum and names such as Dirac and
Sobolev figure prominently.
A good introduction is \cite{friedlander}, while \cite{schwartz} is still
an important work in the field.

Distributions are defined as continuous linear functionals on certain vector
spaces.  Define the space of {\it test functions} by
{${\cal D} =C^\infty_c =\{\phi\fn\R\to\R\mid \phi\in C^\infty \mbox{
 and } \phi \mbox{ has compact support\}}$.  The {\it support} of a function
is the closure of the set on which it does not vanish.  With the usual pointwise
operations $\D$ is a vector space.  An example of a test function
is $\phi(x)=\exp(1/(|x|-1))$ for $|x|<1$ and $\phi(x)=0$, otherwise.  The
only analytic function in ${\cal D}$ is $0$.  We say a sequence $\{\phi_n\}\subset\D$
converges to $\phi\in \D$ if there is a compact set $K$ such that all
$\phi_n$ have support in $K$ and for each integer $m\geq 0$, the sequence
of derivatives $\phi_n^{(m)}$ converges to $\phi^{(m)}$ uniformly on $K$. 
The distributions are then defined as the dual space of $\D$, i.e., the
continuous linear functionals on $\D$. For each $\phi\in\D$, the action
of distribution $T$ is denoted $\langle T,\phi\rangle\in\R$. Linear means that
for all $a,b\in\R$ and all $\phi, \psi\in\D$ we have 
$\langle T,a\phi+b\psi\rangle=a\langle T,\phi\rangle+b\langle T,\psi\rangle$.  
Continuous means  that if
$\phi_n\to\phi$ in $\D$ then 
$\langle T,\phi_n\rangle\to\langle T,\phi\rangle$ in $\R$.
The space of distributions is denoted $\Dp$.

If $f$ is a locally integrable function then $\langle T_f, \phi\rangle=\int_{-\infty}^\infty
f\phi$ defines a distribution since $\phi\in\D$ has compact support, integrals
are linear and dominated convergence or uniform convergence allows us to
take limits under the integral.  Hence,
for $1\leq p\leq \infty$ all the functions in $L^p$ are distributions.
An example of a distribution that is not given by a function is the 
Dirac distribution.  
It is defined
by $\langle\delta,\phi\rangle =\phi(0)$. 

If
$f\in C^1$ and $\phi$ is a test function 
then integration by parts shows that $\int_{-\infty}^\infty f'\phi
=-\int_{-\infty}^\infty f\phi'$.  For all $T\in \Dp$ we can mimic this
behaviour by defining
the derivative via $\langle T',\phi\rangle=-\langle T,\phi'\rangle$.
With this definition, all distributions have derivatives of all orders and
each derivative is a distribution.   For example, 
$\langle\delta',\phi\rangle= -\langle\delta,\phi'\rangle = -\phi'(0)$.  In
electrostatics, $\delta$ models a point charge and $\delta'$ models a dipole.  
If $T$ is a function, we will write its distributional derivative as
$T'$ and its pointwise derivative as $T'(x)$ where $x\in\R$.
From now on, all derivatives will
be distributional derivatives unless stated otherwise.

If $f\in C^0$ then $T_f$ is a distribution.  We can recover its pointwise
value at any point $x\in\R$ by evaluating the limit $\langle T,\phi_n\rangle$
for a sequence $\{\phi_n\}\subset\D$ such that for each $n$, $\phi_n\geq 0$,
$\int_{-\infty}^\infty \phi_n=1$, and the support of $\phi_n$ tends to $\{x\}$
as $n\to\infty$.  Such a sequence is termed a {\it delta sequence}.

The distributional derivative subsumes pointwise and approximate derivatives
and so is very general.  An integration process that inverts it leads to
a very general integral.

\section{The distributional Denjoy integral}\label{defn}

Denote the extended real numbers by $\Rbar=[-\infty,\infty]$.  We define
$C^0(\Rbar)$ to be the continuous functions such that $\lim_{\infty}F$
and $\lim_{-\infty} F$ both exist in $\R$. To be in $C^0(\Rbar)$, $F$
must have real limits at infinity. We can then define $F(\pm\infty)
=\lim_{\pm\infty}F$. Thus, no definition of $F(x)=e^x$ 
at $\pm\infty$ can put $F$ in $C^0(\Rbar)$.  Similarly with
$G(x)=\sin(x)$.  However, $H(x)=\arctan(x)$ is in $C^0(\Rbar)$ if
we define $H(\pm\infty) =\pm\pi/2$.
Define 
$$\Ct=\{F\in C^0(\Rbar)\mid F(-\infty)=0\}.$$
Note that $\Ct$ is a Banach space with the uniform norm $\|F\|_\infty=
\sup_\R |F|=\max_{\Rbar}|F|$.  We now define the space of integrable 
distributions by 
$$
\Alex =\{f\in{\cal D}'\mid f=F' \mbox{ for some }
F\in \Ct\}.$$
A distribution $f$ is integrable if it is the distributional
derivative of a function $F\in\Ct$, i.e., for all $\phi\in\D$ we have
$\langle f,\phi\rangle
=\langle F',\phi\rangle=-\langle F,\phi'\rangle=-\int_{-\infty}^\infty
F\phi'$.  Since $F$ and $\phi'$ are continuous and $\phi'$ has compact
support, this exists as a Riemann integral.  If $f\in\Alex$ then its
integral is defined as $\int_{-\infty}^\infty f=F(\infty)$.  
An obvious
alternative would have been to take $F\in C^0(\Rbar)$ and then
$\int_{-\infty}^\infty f=F(\infty)-F(-\infty)$.  
The function $F$ is a primitive of $f$.

This definition seems to have been first proposed by P.~Mikusi\'nski and
K.~Ostaszewski \cite{pmikusinski1}.  (See also \cite{ostaszewski1},
\cite{ostaszewski2} and 
\cite{pmikusinski2}.)
Without reference to these papers, it was
developed in detail in the plane by D.D.~Ang, K.~Schmidt, L.K.~Vy \cite{ang}
(and repeated in \cite{angvy}).
Several of our results come from this paper.
All of these papers work with the integral in a compact Cartesian interval.

Notice that if $f\in\Alex$
then $f$ has many primitives in $C^0(\Rbar)$, all differing by a constant,
but $f$ has exactly one primitive in $\Ct$.
If $F_1, F_2\in\Ct$ and $F_1'=f, F_2'=f$, then
linearity of the derivative shows $(F_1-F_2)'=0$.  It is known that
the only solutions of
this differential equation are constants \cite[\S2.4]{friedlander}.  
The condition at $-\infty$
now shows $F_1=F_2$. Hence, the integral is unique.

We can define $\int_a^bf=\int_{-\infty}^b f -\int_{-\infty}^a= F(b)-F(a)$
for all $a,b\in\Rbar$, where $F$ is a primitive of $f$.
The integral is then additive: $\int_a^b f+ \int_b^c f=\int_a^c f$.
Also, for open interval $I\subset\R$, define 
{${\cal D}(I) =\{\phi\fn I\to\R\mid \phi\in C^\infty(I) \mbox{
 and } \phi \mbox{ has compact support in }I\}$.  
We then have the distributions on $I$, ${\cal D}'(I)$, being 
the continuous linear functionals on ${\cal D}(I)$.
If $f\in {\cal D}'(I)$ then $f$ is integrable on ${\overline I}$
if there is $F\in C^0({\overline I})$
such that $F'=f$.  Then $\int_a^bf=F(b)-F(a)$.
It is easy to see that these two definitions of $\int_a^b f$ are
equivalent.  For, if $f\in\Dp$ then $f\in\Dp(I)$ since ${\cal D}(I)\subset
{\cal D}$.  If $F\in C^0({\overline \R})$ with $F'=f$ on ${\cal D}$ 
then we also have
$F\in C^0({\overline I})$ so both definitions give $\int_a^b f=F(b)-F(a)$.
If $f\in\Dp(I)$ then in general we cannot extend $f$ to $\Dp$.  For example,
$f$ could be a function with a non-integrable singularity at an endpoint of
$I$.  However, if $f$ is integrable on ${\overline I}$ then we have
$F\in C^0({\overline I})$ such that $F'=f$ on ${\cal D}(I)$. Write $I=(a,b)$.
Define $G=0$ on $[-\infty,a]$, $G=F-F(a)$ on $I$, $G=F(b)-F(a)$ on $[b,\infty]$.
Then $G\in\Ct$ and $G'=f$ on ${\cal D}(I)$.  And, $G(b)-G(a)=F(b)-F(a)$.

Since the derivative is linear, the operations
$\langle af+g,\phi\rangle = a\langle f,\phi\rangle +\langle g,\phi\rangle$
($a\in \R$; $f,g\in\Alex$; $\phi\in\D$) make $\Alex$ into a vector space
and $\int_{-\infty}^\infty af+g =aF(\infty)+G(\infty)$.  
We will use the convention
that when $f, g, f_1$, etc. are in $\Alex$ then we will denote their 
corresponding primitives in $\Ct$ by upper case letters $F, G, F_1$, etc.

Here are some examples that show the extent of applicability of our definition.

\begin{example}
{\rm
1. If $f$ is Riemann integrable on $[a,b]$ then the Riemann integral
$F(x)=\int_a^x f$ is a Lipshitz continuous function and $F'(x)=f(x)$ at all points
of continuity of $f$.  By Lebesgue's characterisation of the Riemann
integral,  $f$ is continuous almost everywhere.  Hence,
$\langle F',\phi\rangle=\int_{-\infty}^\infty
F'\phi=\int_{-\infty}^\infty f\phi$ since changing $F'$ on a set of
measure zero doesn't affect the value of this last integral.  Therefore,
if $F'(x)=f(x)$ almost everywhere then $F'=f$ on $\D$.  The distributional
integral then contains the Riemann integral.  \\

\noindent
2. If $f\in L^1$ then $F(x)=\int_{-\infty}^x f$ defines 
$F\in AC\cap C^0(\Rbar) \subsetneq \Ct$.
Since $F'=f$ almost everywhere, the distributional
integral then contains the Lebesgue integral.  Note that to define
$L^1$ primitives on the real line we have to include the condition
that $F\in  C^0(\Rbar)$  with $F\in AC$.\\

\noindent
3.  If $f$ is Denjoy integrable, then its primitive is an $ACG*$ function
and by the same reasoning as above, the distributional integral contains
the Denjoy integral. 
This integral includes the improper Riemann and Cauchy--Lebesgue extensions
of the Riemann and Lebesgue integrals, respectively.  The function $F'$ given
at the beginning of Section~\ref{arnauddenjoy} has an improper Riemann
integral.  Only the origin is a point of {\it nonabsolute summability},
i.e., over no open interval containing the origin  is $|F'|$ integrable.
However, the Denjoy integral can integrate functions whose set of points of
nonabsolute summability has positive measure, provided it is nowhere dense 
on the real line.  For such functions
it is impossible to define an integral by limits of integrals
over subintervals as is
done with the improper Riemann and Cauchy--Lebesgue processes.  Denjoy
used a transfinite induction process, which he called {\it totalisation},
to define an integral in terms of limits of Lebesgue integrals.  This
was his second solution to the problem of integrating derivatives.  This
integral turned out to be equivalent to the integral defined using $ACG*$
functions. See \cite{bullen} for references to this history.

Denjoy's third solution to the problem of integrating derivatives
was to define an integration process that integrated the approximate
derivative of $ACG$ functions.  Here $ACG$ is yet another complicated
function class of continuous functions that have some differentiability
properties.  In this case, $ACG*\subsetneq ACG\subsetneq C^0$.  
See \cite{celidze}
or \cite{gordon} for details.  The {\it wide} or {\it generalised} 
Denjoy integral
of $f$ is $\int_a^xf=F(x)-F(a)$ where $F\in ACG$ and $D_{ap}F=f$ almost
everywhere.  Using integration by parts for the wide Denjoy integral
\cite[p.~33]{celidze} we can show that if $D_{ap}F=f$ almost
everywhere then $F'=f$ on $\D$.  Hence, the distributional integral
contains the wide Denjoy integral.\\

\noindent
4. Let $F$ be a continuous function such that $F'(x)$ does not
exist for any $x\in\R$.
Then $F'\in\Alex$ and $\int_{a}^b F'=F(b) - F(a)$ for all $a,b\in\R$.
This example shows the following 
difference between Denjoy and distributional integrals.
If $\int_a^b f$ exists as a Denjoy integral then there is a subinterval
$I\subset[a,b]$ such that $|f|$ is integrable over $I$, i.e., $f\in L^1(I)$.  
The corresponding result is false for the distributional integral since
$F$ would have to be $AC$ on $I$ and thus differentiable almost everywhere
in $I$ but $F$ is differentiable nowhere.\\

\noindent
5. Let $F$ be a continuous, increasing, singular function on $[0,1]$,
such as the Cantor--Lebesgue function.  Then $F'(x)=0$ for  almost 
all 
$x\in[0,1]$.  Since $F$ is of bounded variation 
(see Section~\ref{integrationbyparts}), its derivative is
integrable in the Lebesgue sense and $\int_0^x F'(t)\,dt=0$ for all $x\in[0,1]$.
But, $F'\in \Alex$ and the distributional integral is
$\int_0^x F'=F(x)-F(0)$ for all $x\in[0,1]$.\\

\noindent
6. The distributional Denjoy integral is included in the Riemann--Stieltjes
integral since for any function $F$ we have $\int_a^b dF=F(b)-F(a)$.
A valuable feature of the distributional integral is that it confines
itself to the Banach space $\Alex$ so
we can work directly  with the integrand $F'$ rather than have to
deal with the differential $dF$ or its attendant finitely additive measure.
}
\end{example}

We now consider the Banach space structure of $\Alex$.
For $f\in\Alex$,
define  the Alexiewicz norm by  
$\|f\|=\|F\|_\infty=\sup_\R|F|=\max_{\Rbar}|F|$.
\begin{theorem}
With the Alexiewicz norm, $\Alex$ is a Banach space.
\end{theorem}

\bigskip
\noindent
{\bf Proof:} The fact that $\Alex$ is a vector space follows from the
linearity of the derivative, so we will start by proving that $\|\cdot\|$
is a norm.  Let $f,g\in\Alex$.

(i) First, $\|0\|=\|0\|_\infty =0$.  And, if $\|f\|=0$ then $\|F\|_\infty=0$
so $F(x)=0$ for all $x\in\Rbar$.  But then $F'=0$.

(ii) Let $a\in\R$.  Then $(aF)'=a F'$.
Note that this means $\langle (a F)',\phi\rangle=
-\int_{-\infty}^\infty (a F(x))\phi'(x)\,dx
=-a\int_{-\infty}^\infty F(x)\phi'(x)\,dx
=a\langle F',  \phi\rangle$ for all $\phi\in\D$.
We then have $\|a f\|=\|a F\|_\infty=|a|\|F\|_\infty$.

(iii) Since $(F+G)'=F'+G'$ we get $\|f+g\|=\|F+G\|_\infty\leq \|F\|_\infty
+\|G\|_\infty=\|f\|+\|g\|$.

And, $\Alex$ is a normed linear space.  To show it is complete, suppose
$\{f_n\}$ is a Cauchy sequence in $\|\cdot\|$.  Since we have
$\|F_n-F_m\|_\infty=\|f_n-f_m\|$ it follows that $\{F_n\}$ is Cauchy
in $\|\cdot\|_\infty$.  There is $F\in\Ct$ such that
$\|F-F_n\|_\infty\to0$. But then $\|F'-f_n\|=\|F-F_n\|_\infty\to 0$
so $f_n\to F'$ in $\|\cdot\|$.  Since $F\in \Ct$ we have $F'\in\Alex$.\qed

Three equivalent norms are considered in Theorem~\ref{Equivalent norms}.

The definition of the integral shows that $\Alex$ and $\Ct$ are isometrically
isomorphic \cite{ang}.  
They are isomorphic because a  bijection is given by $f\leftrightarrow F$ where $f\in\Alex$
and $F$ is its primitive in $\Ct$.  This mapping is a linear isometry since
for all $F,G\in\Ct$ and all $a\in\R$,
$(aF+G)'=aF'+G'$ and
$\|f\|=\|F\|_\infty$.  This also shows $\Alex$ is separable and that
$L^1$ and the spaces of Denjoy and wide Denjoy integrable functions are
dense in $\Alex$.

\begin{theorem}\label{densityinAlex}
$\Alex$ is separable and $L^1$ and the spaces of Denjoy and wide Denjoy integrable functions are
dense in $\Alex$.
\end{theorem}

\bigskip
\noindent
{\bf Proof:} 
Functions, $\Phi$, for which there is a polynomial
$p$ and an
interval $[a,b]\subset\R$ with $p(a)=0$ such that $\Phi=0$ on $(-\infty,a]$,
$\Phi=p$ on $[a,b]$, and $\Phi=p(b)$ on $[b,\infty)$, are dense in $\Ct$ with
$\|\cdot\|_\infty$.  But such functions are absolutely continuous, so
$L^1$ is
dense in $\Alex$.  It follows that the spaces of Denjoy and wide Denjoy 
integrable functions are
dense in $\Alex$.  Polynomials on $[a,b]$ with rational coefficients 
form a countable
dense set in $C^0([a,b])$ so $\Alex$ is separable.\qed

Note that $C^0(K)$ is separable exactly when $K$ is compact 
\cite[Exercise~V.7 17]{dunford}.
Our two-point compactification of the real line makes $\Rbar$ into a
compact Hausdorff space.  A topological base is the set of all intervals
$(a,b)$, $[-\infty,b)$, $(a,\infty]$ and $[-\infty,\infty]$, for all
$-\infty\leq a< b\leq \infty$.  That is, we declare all such intervals
open in $\Rbar$.
 
One half of the Fundamental Theorem
is built into the definition.  The other half follows easily.

\begin{theorem}[Fundamental Theorem of Calculus]\label{fundamental}
{\rm (a)} Let $f\in \Alex$. Define $\Phi(x)=\int_{-\infty}^x f$.  Then
$\Phi\in \Ct$ and $\Phi'=f$.\\
{\rm (b)}  Let $F\in C^0(\Rbar)$.  Then $\int_{-\infty}^x F'=F(x)-F(-\infty)$ for all
$x\in \Rbar$.
\end{theorem}

\bigskip
\noindent
{\bf Proof:} (a) By the uniqueness of the integral, $\Phi=F\in \Ct$.  Then
$\Phi'=F'=f$.  (b) There is a constant $c\in\R$ such that $F+c\in\Ct$.  But then
$F'=(F+c)'\in\Alex$ and the result follows from the definition of the integral.
\qed

At this stage it is hoped that the reader appreciates what we have 
accomplished.  With minimal effort we have proven a very general version of the 
Fundamental Theorem and have proven that the space of integrable 
distributions is a Banach space.  
To prove the corresponding results for the Lebesgue integral requires
considerably more machinery.  For example, part~(b) of 
Theorem~\ref{fundamental} (Lebesgue differentiation theorem) uses
the Vitali covering theorem.  And, one usually requires convergence
theorems to prove that $L^1$ is complete.  

\section{Integration by parts,  H\"older's inequality}\label{integrationbyparts}

If $g\fn\Rbar\to\R$, its {\it variation} is $Vg=\sup\sum_n|g(x_n)-g(y_n)|$
where the supremum is taken over every sequence $\{(x_n,y_n)\}$ of 
disjoint intervals in $\Rbar$.  The set of functions with bounded variation
is denoted  $\bv$.
It is known that functions of bounded variation
are bounded and
have left and right limits at each point (from the right at $-\infty$
and from the left at $\infty$.)  Thus, if $g\fn\R\to\R$ is of bounded
variation on $\R$ then the limits $\lim_{-\infty}g$ and $\lim_\infty g$
exist and we will use these to extend the domain of $g$ to $\Rbar$.  If $g\in\bv$ we can change $g$ on a
countable set so that it is right continuous on $[-\infty,\infty)$ and
left continuous at $\infty$, i.e., 
$\lim_{x\to a^+}g(x)=g(a)$ for all $a\in[-\infty,\infty)$
and 
$\lim_{x\to\infty}g(x)= g(\infty)$.
We will say such  functions are of {\it normalised bounded variation} 
(${\cal NBV}$).  (This is slightly different from the usual definition
but more convenient for our purposes.  See \cite[p.~241]{dunford}.)  
The space $\bv$ is a Banach space with norm
$\|g\|_{\bv}=|g(-\infty)|+Vg$.

The {\it essential variation} is defined as
$essvar\, g=\sup\int_{-\infty}^\infty g\phi'$ where now the supremum
is taken over all $\phi\in C^1_c$ with $\|\phi\|_\infty\leq 1$.
Denote the functions of essential variation by $\ebv$.
Changing a function at even one point can affect its variation but
changing a function on a set of measure zero will not affect its
essential variation.  And, $\bv\subsetneq\ebv$ but changing a function
in $\ebv$ on a certain set of measure zero will put it into $\bv$.
The space $\ebv$ is a Banach space with norm
$\|g\|_{\ebv}=\|g\|_{\infty}+essvar \,g$.  
If $g\in{\cal NBV}$ then its variation and essential variation are 
identical.  

As with the Denjoy integral, 
functions of bounded
variation play an important role in the distributional integral.
They form the dual space, tell us about integration by parts and
H\"older's inequality. Theorem~\ref{measurezero} below shows that
results that hold for functions of bounded variation also hold
for functions of essential bounded variation.

If $F\in C^0(\Rbar)$ and $g\in\bv$ then it is known that the
Riemann--Stieltjes integral $\int_{-\infty}^\infty F\,dg$ exists.  It
can be defined using a
partition of $\Rbar$.
The integral exists, with value $\int_{-\infty}^\infty F\,dg\in\R$,
if for all $\epsilon >0$ there is $\delta>0$ 
such that if
$-\infty=x_0<x_1<\ldots
<x_N=\infty$, 
$\max_{2\leq n\leq N-1}(x_n-x_{n-1})<\delta$,
$x_1<-1/\delta$, and  $x_{N-1}> 1/\delta$  then, for all $z_n\in[x_{n-1},x_n]$,
we have 
$$
\left|\sum_{n=1}^N F(z_n)\left[g(x_n)-g(x_{n-1})\right] -\int_{-\infty}
^\infty F\,dg\right|<\epsilon.
$$
To integrate over $[a,b]\subset\Rbar$ we use partitions of $[a,b]$.
The integral can also be defined by taking limits of 
Riemann--Stieltjes integrals
over finite subintervals:
\begin{eqnarray*}
\int_{-\infty}^\infty F\,dg & = & \lim_{\stackrel{B\to\infty}{A\to-\infty}}
\int_A^BF\,dg+F(\infty)\left[g(\infty)-\lim_{x\to\infty}g(x)\right]\\
 & & \quad + F(-\infty)\left[\lim_{x\to-\infty}g(x)-g(-\infty)\right].
\end{eqnarray*}
See \cite[p.~187]{mcleod} and \cite{talvilarae} for details.  

\begin{prop}\label{propH}
Let $f\in\Alex$ and $g\in\bv$.  Define
$H(x)=F(x)g(x)-\int_{-\infty}^x
F(t)\,dg(t)$.  Then $H\in\Ct$.
\end{prop}

\bigskip
\noindent
{\bf Proof:}  Since $g$ is of bounded variation, it is bounded.  Write
$|g|\leq M$ for some $M\in\R$.
Let $x\in\R$ and $y\geq x$.    Because $\int_x^y dg=g(y)-g(x)$,
we can write
$
H(x)-H(y)
 =  \left[F(x)-F(y)\right] g(x) + \int_{x}^y\left[F(t)-F(y)\right]dg(t)
$.
Now,
\begin{eqnarray}
|H(x)-H(y)|  
 & \leq & |F(x)-F(y)|M + \max_{x\leq t\leq y}|F(t)-F(y)|Vg\label{ibyp}\\
 & \to & 0 \mbox{ as } y\to x \mbox{ since } F \mbox{ is uniformly continuous.}
\notag
\end{eqnarray}
Similarly if $y\leq x$.  Hence, $H\in C^0(\R)$.  To prove $H\in\Ct$, let
$x\in\R$.  Then $|H(x)|\leq |F(x)|M+\|F\chi_{(-\infty,x]}\|_\infty
 Vg\to 0$ as $x\to-\infty$.  From \eqref{ibyp}, the sequence
$\{H(n)\}$ is Cauchy and so has a limit as $n\to\infty$.
Hence, $H\in\Ct$. \qed

We now get the integration
by parts formula.
\begin{defn}[Integration by parts]
Let $f\in\Alex$ and $g\in\bv$.  Define $fg=H'$ where
$H(x)=F(x)g(x)-\int_{-\infty}^x F\,dg$.  Then $fg\in\Alex$ and
 $\int_{-\infty}^\infty fg=F(\infty)g(\infty) -
\int_{-\infty}^\infty F\,dg$.
\end{defn}
Notice that in 
$H(x)=F(x)g(x)-\int_{-\infty}^xF\,dg$ we really mean $g(x)$ and
not the left or right limit of $g$ at $x$, 
including the cases when $x=\pm\infty$.  
Although $g$ has a limit at infinity, it
might also have a jump discontinuity at infinity.  
Changing $g$ on a countable set will in general change the value
of both $F(x)g(x)$ and $\int_{-\infty}^xF\,dg$ but will not affect
$H(x)$.  To see this, it suffices to prove that if $g\in\bv$ and
$g=0$, except perhaps on a countable set,  then $H=0$.  Let $x\in\R$ and $\epsilon>0$.  Since
$\lim_{-\infty}F=0$, we can take
$A<x$ such that  $g(A)=0$ and
$|\int_{-\infty}^A F\,dg|\leq\|F\chi_{(-\infty,A]}\|_{\infty}Vg<\epsilon/3$.
Since $F$ is continuous at $x$, we can take $A<B<x$ such that $g(B)=0$ and
\begin{eqnarray*}
\left|F(x)g(x)-\int_{B}^x F\,dg\right| & = & 
\left|F(x)[g(x)-g(B)] -\int_{B}^x F\,dg\right|\\
 & = & \left|\int_B^x\left[F(x)-F(t)\right]dg(t)\right|\\
 & \leq & \max_{B\leq t\leq x}|F(x)-F(t)|\,Vg\\
 & \leq & \epsilon/3.
\end{eqnarray*}
And, since $F$ is uniformly continuous, there are $A=a_0<a_1<\ldots<
a_N=B$ such that $g(a_n)=0$ for all $0\leq n\leq N$ and
$\max_{a_{n-1}\leq t\leq a_n}|F(a_n)-F(t)|<\epsilon/[3(1+Vg)]$.  Then
\begin{eqnarray*}
\left|\int_A^BF\,dg\right| & = & \left|\sum_{n=1}^N\int_{a_{n-1}}^{a_n}
F\,dg\right|\\
 & = & \left|\sum_{n=1}^N\int_{a_{n-1}}^{a_n}\left[F(a_n)-F(t)\right]
dg(t)\right|\\
 & \leq & \sum_{n=1}^N\max_{a_{n-1}\leq t\leq a_n}|F(a_n)-F(t)|\,V(g\chi
_{[a_{n-1},a_n]})\\
 & \leq & \frac{\epsilon\,Vg}{3(1+Vg)}.
\end{eqnarray*}
Combining these results shows that $H(x)=0$.

A general distribution $T\in {\cal D}'$ can be multiplied by a smooth
function $h\in C^\infty$ using
$\langle hT,\phi\rangle=\langle T, h\phi\rangle$.
This works because $h\phi\in {\cal D}$ for all $\phi\in{\cal D}$.  
We can multiply
$f\in \Alex$ 
by any function $g\in\bv$.  Define
$fg=H'$, i.e.,
\begin{eqnarray*}
\langle fg, \phi\rangle & = & \langle H', \phi\rangle  = 
-\langle H, \phi'\rangle\\
 & = & -\int_{-\infty}^\infty\left[F(x)g(x)-\int_{-\infty}^x F(t)\,dg(t)
\right]\phi'(x)\,dx.
\end{eqnarray*}
Since $\phi$ is of compact support, Fubini's theorem tells us we can
interchange orders of integration to write
$\langle fg, \phi\rangle=\langle (Fg)',\phi\rangle
-\int_{-\infty}^\infty F(t)\phi(t)\,dg(t)$.  This agrees with the
usual definition when $g\in C^\infty$ since then for
$\phi\in\D$ we have $g\phi\in\bv$ and  $\langle fg,\phi\rangle=
\langle f,g\phi\rangle =\langle (Fg)',\phi\rangle
-\int_{-\infty}^\infty F(t)\phi(t)\,dg(t)$, upon integrating by parts.

The integration by parts formula agrees with the usual one when
$f$ has a Lebesgue, Henstock--Kurzweil or wide Denjoy integral.  
Note that we 
have defined $fg=H'$ but we have no way of proving this.  However, we can
use the norm to show this is the correct definition.  Suppose $f\in\Alex$
and $g\in\bv$ with $|g|\leq M$. 
By Theorem~\ref{densityinAlex}, 
there is a sequence $\{f_n\}\subset L^1$ such that
$\|f_n-f\|\to 0$ as $n\to\infty$.  Define $H_n(x):=\int_{-\infty}^x
f_ng=F_n(x)g(x)-\int_{-\infty}^x F_n\,dg$ by the usual integration by
parts formula.  As in \eqref{ibyp},
$|H_n(x) -H(x)| \leq \|F_n-F\|_\infty(M+Vg)\to 0$ as $n\to\infty$.
It follows that $\|H_n-H\|_\infty\to 0$, which justifies our definition
$fg=H'$.

Note that for $(a,b)\subset\R$ we have 
$\int_a^b fg=F(b)g(b)-F(a)g(a)-\int_a^bF\,dg$,
where $F'=f$ and $F\in C^0([a,b])$.  A consequence is that if $f\in\Alex$
then $f$ is integrable on every subinterval of the real line.  For
compact interval $[a,b]$,
\begin{eqnarray*}
\int_a^b f & = & \int_{-\infty}^\infty f\chi_{[a,b]}\\
 & = & F(\infty)\chi_{[a,b]}(\infty)-
\int_{-\infty}^\infty
F\,d\chi_{[a,b]}\\
& = & F(b)-F(a).
\end{eqnarray*}
We also have $\int_{I}f=
F(b)-F(a)$ when $I=[a,b]$, $[a,b)$, $(a,b]$ or $(a,b)$.  This can be seen by
letting $g=\chi_{I}$ and integrating by parts.  

The integration by parts formula shows that the distributional integral is compatible with Schwartz's definition of integral \cite{schwartz}.  If
$f\in\Dp$ such that $f(1)$ is defined then $\int_{-\infty}^\infty f :=f(1)$.  
Since the function $1\in\bv$, integration by parts gives,
$f(1)=\int_{-\infty}^\infty f=\int_{-\infty}^\infty f1=F(\infty)1-
\intinf f\,d1=F(\infty)$.  For another type of distributional integral, see
the final paragraph of
Section~\ref{odds}.

As a corollary to Proposition~\ref{propH} we have a version of the
H\"older inequality.  

\begin{theorem}[H\"older inequality]\label{holder}
Let $f\in\Alex$. If $g\in{\cal NBV}$ then
$\left|\intinf fg\right|  \leq  |\intinf f|\inf_{\R}|g|+2\|f\|Vg$.
 If $g\in{\cal BV}$ then
$\left|\intinf fg\right|  \leq  2\|f\| \|g\|_{\bv}$.
\end{theorem}
The first inequality was proved in \cite[Lemma~24]{talvilafourier} for the
Henstock--Kurzweil integral and the same proof works here.  
The second inequality  is similar.  The factor of `2' is replaced by
`1' if we use the equivalent norm on $\Alex$, $\|f\|':=\sup_I|\int_I f|$
where the supremum is taken over all intervals $I\subset\R$.

We now get a new interpretation of the action of $f\in\Alex$
as a distribution.  Let $\phi\in\D$.  Since $\phi\in\bv\subset C^1$, we have
\begin{eqnarray*}
\langle f,\phi\rangle & = & \langle F',\phi\rangle=-\langle F,\phi'\rangle
=-\intinf F\phi'\\
 & = & -\intinf F\,d\phi = \intinf f\phi -F(\infty) \phi(\infty) \\
 & = & \intinf f\phi.
\end{eqnarray*}
Hence, the action of $f$ on test function $\phi$ is interpreted as
the integral of the product $f\phi$, as in the case when $f$ is a 
locally integrable function.

The H\"older inequality shows that $f$ is a continuous linear functional
on $\bv$.  Suppose $\{g_n\}\subset \bv$ and $\|g_n\|_\bv\to 0$ as
$n\to\infty$.  Then $f$ is continuous:
$$
|\langle f,g_n\rangle| = \left|\intinf fg_n\right|\leq 
2\|f\| \|g_n\|_{{\cal BV}}\to 0.
$$  
And, for $a\in\R$; $g_1, g_2\in\bv$; 
\begin{eqnarray*}
\langle f,ag_1+g_2\rangle & = & F(\infty)\left[ag_1+g_2\right]\!(\infty)
-\intinf F\,d(ag_1+g_2)\\
 & = & aF(\infty)g_1(\infty) + F(\infty)g_2(\infty)-a\intinf F\,dg_1
-\intinf F\,dg_2\\
 & = & a\langle f,g_1\rangle +\langle f,g_2\rangle.
\end{eqnarray*}
So, we know that the dual of $\bv$ contains $\Alex$, i.e., 
$\Alex\subset\bv^*$.  In fact,
$\bv^*$ is much larger than $\Alex$ since it contains measures not
in $\Alex$ such as the Dirac measure.
However, we do know that $\Alex^*=\bv$.  If
$\{f_n\}\subset\Alex$ and $\|f_n\|\to 0$ then for $g\in\bv$ it follows that
$\left|\intinf f_n g\right| \leq 2\|f_n\|\|g\|_\bv \to 0$ so $g\in\Alex^*$,
since we also have linearity $\langle af_1+f_2, g\rangle=
a\langle f_1,g\rangle +\langle f_2,g\rangle$.
The Riesz Representation Theorem says that if $[a,b]$ is a compact interval then
$C^0([a,b])^*=\bv$.  Since our two-point compactification of the real
line makes $\Ct$ homeomorphic to the continuous functions on $[a,b]$
vanishing at $a$, it also true that $\Alex^*=\bv$.
Hence, the functions of bounded variation are the multipliers
for the distributional integral ($g\in\bv$ implies $fg\in\Alex$ for all
$f\in\Alex$) and $\bv$ also forms the dual space (the set of continuous linear functionals
on $\Alex$).  

Although it is prohibited
to discuss measure and distribution $f\in\Alex$ in the same breath, 
measure-theoretic arguments apply to $g\in\bv$. Using a density 
argument, we see that changing $g$ on a set of
measure 0 does not affect the value of $\intinf fg$.
\begin{theorem}\label{measurezero}
Let $f\in\Alex$ and let $g\in\cal{EBV}$.
Let
$\{\phi_n\}\subset\D$ with $\|f-\phi_n\|\to 0$. 
Define $\intinf fg=\lim_{n\to\infty}\intinf \phi_ng$.
Let $\tilde g$ be the unique function in $\cal{NBV}$ such that $essvar \,g
=V
{\tilde g}$.
Then $\intinf fg=\intinf f\tilde{g}$.
\end{theorem}

\bigskip
\noindent
{\bf Proof:}
Note that such a sequence $\{\phi_n\}$ exists since $\D$ is dense in
$\Alex$.  For each $n\in\N$, the integral $\intinf \phi_ng$ exists
as a Lebesgue integral
since $\phi_n$ is smooth with compact support and $g\in L^1_{loc}$.
We can then change $g$ on a set of measure zero to get
$\intinf \phi_ng=\intinf \phi_n\tilde{g}\to \intinf f\tilde{g}$, using a convergence
theorem for Henstock--Kurzweil integrals \cite[Corollary~3.3]{talvilarae}.
The definition does not depend on the choice of $\{\phi_n\}$ since
if $\{\psi_n\}\subset\D$ with $\|f-\psi_n\|\to 0$ then
\begin{eqnarray*}
\left|\intinf \phi_ng-\intinf \psi_ng\right| & = & \left|\intinf(\phi_n-\psi_n)\tilde{g}\right|\\
 & \leq & 2\|\phi_n-\psi_n\|\|{\tilde g}\|_{\bv}\to 0\quad\text{as } n\to\infty.
\end{eqnarray*}
Hence we are justified in writing $\intinf fg=\intinf f\tilde{g}$ for all $f\in
\Alex$.\qed

\begin{corollary}
$\Alex^*=\ebv$.
\end{corollary}

The H\"older inequality also shows that if $f\in\Alex$ then $f$ is
a distribution of order one and hence is tempered.  See \cite{friedlander} for
the definitions.

\section{Change of variables}

In order to write a change of variables formula, we need to be able
to compose a distribution in $\Alex$ with a function.  For $(\alpha,\beta)
\subset\R$, we can define
$\D((\alpha,\beta))$ to be the test functions with compact support in $(\alpha,
\beta)$ and then $\Dp((\alpha,\beta))$ is the corresponding space of 
distributions.  Suppose $(\alpha,\beta), (a,b)\subset\R$. 
If we have distribution $T\in\Dp((\alpha,\beta))$, let 
$G\fn(a,b)\to(\alpha,\beta)$ be a $C^\infty$ bijection such that
$G'(x)\not=0$ for any $x\in(a,b)$.  Then $T\circ G\in\D'((a,b))$ is
defined by $\langle T\circ G,\phi\rangle =
\langle T,\frac{\phi\circ G^{-1}}{G'\circ G^{-1}}\rangle$ for all
$\phi\in\D((a,b))$.  This definition follows from the change of variables
formula for smooth functions.  See \cite[\S7.1]{friedlander}.
For $f\in\Alex$ and $G$ as above, this
then leads to the formula
$
\int_\alpha^\beta f=\int_a^b (f\!\circ\! G)\, G'
$ when $G$ is increasing, with a sign change if $G$ is decreasing.
However, using the properties of $\Alex$, we can do much better than this.  
We will show below that the norm validates this formula 
when the
only condition on $G$ is that it be continuous.  First we need to
define the derivative of the composition of two continuous functions.

\begin{defn}[Derivative of composition of continuous functions]
\label{composition}
Let $F,G\in C^0(\Rbar)$.  Then $(F'\circ G)G'  := (F\circ G)'$, i.e.,
$\langle (F'\circ G)G',\phi\rangle =\langle(F\circ G)',\phi\rangle
=-\langle F\circ G,\phi'\rangle = -\intinf (F\circ G)(t)\,\phi'(t)\,dt
$ for all $\phi\in\D$.
\end{defn}

The Alexiewicz norm shows this definition is compatible with the usual
definition for smooth functions. Suppose $F,G\in C^0(\Rbar)$. 
Let $\epsilon>0$.  Take $\delta>0$
such that whenever $|x-y|<\delta$ we have $|F(x)-F(y)|<\epsilon/2$.
This is possible since $F$ is uniformly continuous on $\Rbar$.
There are
$C^1$ functions $p$ and $q$ such that
$\|F-p\|_\infty<\epsilon/2$ and 
$\|G-q\|_\infty<\delta$.  Note that
$F\circ G\in C^0(\Rbar)$ so $(F\circ G)'\in\Alex$.  And,
$(p\circ q)'(t)=(p'\circ q)(t)\,q'(t)$ for all $t\in\R$.  We have
\begin{eqnarray*}
\|(F\circ G)'-(p'\circ q)q'\| & = & \|F\circ G-p\circ q\|_\infty\\
 & \leq & \|F\circ G -F\circ q\|_\infty + \|(F-p)\circ q\|_\infty\\
 & < & \epsilon/2 + \epsilon/2.
\end{eqnarray*}

With this definition we then have the following change of variables
formula.

\begin{theorem}\label{changeofvariablestheorem}
Suppose $f\in\Alex$ and $F'=f$ where $F\in C^0(\Rbar)$.  Let
$-\infty\leq a<b\leq\infty$. If
$G\in C^0([a,b])$  then 
$$
\int_{G(a)}^{G(b)} f  =  \int_a^b (f\circ G)\,G'=
(F\circ G)(b) -(F\circ G)(a).
$$
If $G\in C^0((a,b))$ and $\lim_{t\to a^+}G(t)=-\infty$ and
$\lim_{t\to b^-}G(t)=\infty$ then 
$$
\intinf f = \int_a^b(f\circ G)\,G'
=F(\infty)-F(-\infty).
$$
\end{theorem}

The first statement follows from Definition~\ref{composition} and the
second from Theorem~\ref{haketheorem} below.
This is a remarkable formula because it demands so little of $f$ and
$G$.  For Lebesgue integrals, the usual formula requires $f\in L^1$ and
$G\in AC$ and monotonic \cite[\S38.4]{mcshane}.  
Even invoking Stieltjes  integrals leads to a change of variables
formula requiring monotonicity or differentiability properties of $G$.
See \cite[Exercises~III.13 4. and 5.]{dunford}.  Similarly for the 
Denjoy integral.  See \cite[\S2.7, \S7.9]{mcleod}.  J.~Foran \cite{foran} cites 
references to further theorems in Denjoy integration.  See \cite{bagby}
and \cite{sarkhel}
for good change of variables theorems for Riemann integrals.

\section{Convergence Theorems}\label{sectionconvergence}

Two of the main reasons the Lebesgue integral so easily replaced the
Riemann integral in the first part of the twentieth century were that
the space $L^1$ is a Banach space and there are excellent convergence
theorems.  We have already shown that $\Alex$ is a Banach space.  Now
we will look at convergence theorems.

A sequence $\{f_n\}\in \Alex$ is said to converge {\it strongly} to
$f\in\Alex$ if $\|f_n-f\|\to 0$.  It converges {\it weakly} in
$\D$ if $\langle f_n-f,\phi\rangle =\intinf (f_n-f)\phi\to 0$ for each
$\phi\in\D$.  And, $\{f_n\}$ converges   {\it weakly} in
$\bv$ if $\intinf (f_n -f) g\to 0$ for each
$g\in\bv$. 
\begin{theorem}
Weak convergence in $\bv$ implies weak convergence in $\D$.
Strong convergence implies weak convergence in $\D$ and $\bv$.
Weak convergence in $\D$ does not imply weak convergence in $\bv$.
Weak convergence in $\bv$ does not imply strong convergence.
\end{theorem}

\bigskip
\noindent
{\bf Proof:} Since $\D\subset\bv$, weak convergence in $\bv$ implies 
weak convergence in $\D$.  Suppose $\|f_n-f\|\to 0$.  Then
$\|F_n-F\|_\infty\to 0$.  Let $g\in\bv$.  By the H\"older inequality, 
$$
\left|\langle f_n-f ,g\rangle\right| = \left|\intinf(f_n-f)\,g\right|
\leq 2\|F_n-F\|_\infty \|g\|_{\bv}\to 0.
$$
To see that weak convergence in $\D$ does not imply weak convergence
in $\bv$, let $f_n=\chi_{(n, n+1)}$. For $\phi\in\D$ we have
$\intinf f_n\phi = \int_n^{n+1}\phi \to 0$ but if $g=1$ then 
$\intinf f_ng =1\not\to 0$.
Let $f_n=\chi_{(n-1,n)}-\chi_{(n,n+1)}$.
Then $\{f_n\}\subset\Alex$.  For $g\in \bv$, we have
$\intinf f_ng\to 0$ by dominated convergence
since
$\|f_n\|_\infty=1$, $g$ is bounded and $f_n\to 0$ pointwise on $\R$.
As $\|f_n\|=1$, weak convergence
in $\bv$ (and hence in $\D$) does not imply strong convergence.\qed

Suppose $\{f_n\}\subset\Alex$.  
Strong convergence $\|f_n-f\|\to 0$ implies $f\in\Alex$
since $\Alex$ is a Banach space.  If $f_n\to f$ weakly in $\bv$
then by definition $f\in\Alex$.  But, if $f_n\to f$ weakly
in $\D$ then $f$ need not be in $\Alex$.

\begin{example}
{\rm
There is a sequence $\{f_n\}\subset
\Alex$ that converges weakly in $\D$ to $f\in\Dp\setminus\Alex$.
Let $f_n=\chi_{[-n,n]}$.  Then $f_n\in\Alex$ for each $n\in\N$.
Let $\phi\in\D$.  By dominated convergence (or Weierstrass $M$-test),
$\langle f_n,\phi\rangle
=\int_{{\rm supp(\phi)}}\chi_{[-n,n]}\phi\to\intinf \phi=\langle 1,\phi\rangle$.
Hence, $f_n$ converges weakly in $\D$ to $1\in\Dp\setminus\Alex$.\qed
}
\end{example}

Now suppose we are interested in conditions on $f_n$ so that
$\intinf f_n\to \intinf f$.  

\begin{theorem}\label{theoremconv2}
Let $\{f_n\}\subset\Alex$ and $f\in\Alex$.
If $\|f_n-f\|\to 0$ then $\intinf f_n\to \intinf f$.  The converse is false.
If $f_n\to f$ weakly in $\bv$ then $\intinf f_n\to \intinf f$.  There is a
sequence $\{f_n\}\subset\Alex$ and a distribution $f\in\Alex$ such that $f_n\to f$ weakly in $\D$ 
and $\intinf f_n \not \to \intinf f$.
There is a sequence $\{f_n\}\subset \Alex$ that does not converge weakly
in $\D$ but $\{\intinf f_n\}$ converges in $\R$.
\end{theorem}

\bigskip
\noindent
{\bf Proof:} Certainly we have $|\intinf f_n-f|\leq \|F_n-F\|_\infty=\|f_n-f\|$
so $\|f_n-f\|\to 0$ and the triangle inequality imply 
$\intinf f_n\to \intinf f$.  
Let $f_n(t)=n^2\sin(nt)$ for $|t|\leq \pi$ and $f_n(t)=0$ for $|t|>\pi$.
Then for each $n\in\N$, $\intinf f_n=0$ but $\|f_n\|=n^2\int_{0}^{\pi/n}
\sin(nt)\,dt=2n\to\infty$.  Now suppose $f_n\to f$ weakly in $\bv$.  Since
$1\in\bv$ we have $\intinf f_n\to \intinf f$.  And, define
$$
F_n(t)=\left\{\begin{array}{cl}
0, & t\leq n\\
t-n, & n\leq t\leq n+1\\
1, & t\geq n+1.
\end{array}
\right.
$$
Then $F_n\in\Ct$ and 
$
f_n(t):=F_n'(t)= 1$ for $
n\leq t\leq n+1$ and $f_n(t) =0$, otherwise.
For $\phi\in {\cal D}$, 
$
\langle f_n,\phi\rangle = \int_{n}^{n+1}\phi\to 0$
since $\phi$ has compact support.  But, $\intinf f_n=F_n(\infty)=1$.
This phenomenon can also occur on compact intervals.  Let $F_n(t)=t^n$ for
$t\in[0,1]$.  Then $\int_0^1 f_n=F_n(1)=1$ and yet, for $\phi\in\D((0,1))$,
$|\langle f_n,\phi\rangle|=|-\int_0^1t^n\phi'(t)\,dt|\leq\|\phi'\|_\infty
\int_0^1 t^n\,dt
=\frac{\|\phi'\|_\infty}{n+1}\to 0$.
Finally, let  $f_n(t)=a_n$ for $1\leq t\leq 2$,
$f_n(t)=-a_n$ for $-2\leq t\leq -1$, and $f_n(t)=0$, otherwise.  Here,
$\{a_n\}$ is an arbitrary sequence of real numbers.  Then, $\intinf f_n=0$
for each $n\in\N$ but, unless $\lim_{n\to\infty}a_n=0$, $\{f_n\}$ is
not weakly convergent in $\D$ since we can always take a test function
that has support in $[0,3]$ that is identically 1 on $[1,2]$. \qed

Theorem~\ref{theoremconv2} indicates that to have $\intinf f_n\to \intinf f$ we
should look for some condition between weak convergence in $\bv$, which
is sufficient but not necessary, 
and weak convergence in $\D$, which is neither necessary nor sufficient.  Note
that for $\intinf f_n\to \intinf f$ we will really want $F_n(x)\to F(x)$ 
for each $x
\in\Rbar$.  Indeed,  a corollary to Theorem~\ref{theoremconv2} is that
strong convergence or weak convergence  in $\bv$ of $f_n\to f$ both imply
$\int_{-\infty}^x f_n\to \int_{-\infty}^x f$ for all $x\in\Rbar$.
If we do not have convergence on subintervals then
each $f_n$ could be an arbitrary distribution in $\Alex$ with integral 0 
and we would
then not expect there to be any sensible condition on $f_n$ that ensures
$\intinf f_n\to 0$.

Note that strong convergence $\|f_n-f\|\to 0$ is the same as uniform
convergence of $F_n\to F$ on $\Rbar$.  If each function $F_n\in\Ct$ then uniform
convergence of $F_n\to F$ guarantees $F$ is continuous on $\Rbar$.
Since each $F_n(-\infty)=0$, we also have $F(-\infty)=0$ so $F\in\Ct$
and $\int_{-\infty}^x f_n\to\int_{-\infty}^x F'$ for each $x\in\Rbar$.
But, uniform convergence is not necessary for the limit of a sequence of
continuous functions to be continuous.  The necessary and sufficient
condition is {\it quasi-uniform} convergence.  See \cite{gordon2} or
\cite[IV.6.10]{dunford}.
\begin{defn}[Quasi-uniform convergence]
Let $\{F_n\}\subset C^0(\Rbar)$ and suppose $F\fn\Rbar\to\R$.
If $F_n(x)\to F(x)$ at
each point $x\in\Rbar$ then $F_n\to F$ quasi-uniformly at $x\in\R$ if for each
$\epsilon >0$ and each $N\in\N$ there is $\delta>0$ and $n\geq N$ such that
whenever $|x-y|<\delta$ we have $|F_n(y)-F(y)|<\epsilon$.  For quasi-uniform
convergence at $x=\infty$, replace the condition involving $\delta$ with
$y>1/\delta$, with a similar condition for $x=-\infty$.
\end{defn}

\begin{theorem}\label{quasiuniform}
Let $\{f_n\}\subset\Alex$ and $F\fn \Rbar\to
\R$.
If $F_n\to F$ quasi-uniformly on $\Rbar$ then $F\in\Ct$ and
$\int_{-\infty}^x f_n\to \int_{-\infty}^x F'$ for each $x\in\Rbar$.
\end{theorem}

The following three results give sufficient conditions for 
$\int_{-\infty}^x f_n$ to converge to $\int_{-\infty}^x f$.
Each involves weak convergence of $f_n\to f$ in $\D$.

\begin{theorem}[\cite{ang}, Theorem~8]\label{ang8}
Let $\{f_n\}\subset\Alex$ and
$F\in C^0(\Rbar)$.  Suppose $\{F_n\}$ is uniformly bounded
on each compact interval in $\R$ and $F_n\to F$ on $\Rbar$.  Then
$f_n\to F'$ weakly in $\D$ and $\int_{-\infty}^x f_n\to\int_{-\infty}^x F'$
for each $x\in \R$.
\end{theorem}

\bigskip
\noindent
{\bf Proof:} Since $F(-\infty)=\lim_{n\to\infty}F_n(-\infty)=\lim_{n\to\infty}
0=0$ we have $F\in\Ct$.  Let $\phi\in\D$ with support in the compact
interval $I\subset\R$.  Then $|\langle F_n,\phi\rangle|= |\int_I F_n\phi|
\leq \|F_n\phi\chi_I\|_\infty\lambda(I)$.  By dominated convergence (or the
Weierstrass $M$-test), $\intinf F_n\phi \to \intinf F\phi$, i.e., $F_n\to F$
weakly in $\D$.  And, since $\phi'\in\D$,
$\langle f_n,\phi\rangle=-\langle F_n,\phi'\rangle\to -\langle F,\phi'\rangle
=\langle F',\phi\rangle$.  Therefore, $f_n\to F'$ weakly in $\D$.  And,
$\int_{-\infty}^x f_n=F_n(x)\to F(x) =\int_{-\infty}^x F'$ for each
$x\in\Rbar$.\qed

\begin{corollary}[\cite{ang}, Theorem~9]\label{ang9}
Let $\{f_n\}\subset\Alex$ and
$F\in C^0(\Rbar)$.  Suppose $\{F_n\}$ is uniformly bounded
on each compact interval in $\R$ and $F_n\to F$ on $\Rbar$.  
Suppose $f_n\to f$ weakly in $\D$ for some $f\in \Dp$.  Then $f=F'\in\Alex$
and
$\int_{-\infty}^x f_n\to\int_{-\infty}^x f$
for each $x\in \Rbar$.
\end{corollary}

\bigskip
\noindent
{\bf Proof:} As in the theorem, $F_n\to F$ weakly in $\D$. Therefore,
for $\phi\in\D$,
$\langle f_n,\phi\rangle=-\langle F_n,\phi'\rangle\to-\langle F,\phi'\rangle$.
By the uniqueness of limits in $\D$, $f=F'\in\Alex$. \qed

A sequence of functions $\{F_n\}\subset
\Ct$  is {\it equicontinuous} at $x\in\R$ if for all $\epsilon>0$ there exists
$\delta>0$ such that for all $n\geq 1$, if $y\in\R$ such that
$|x-y|<\delta$ then $|F_n(x)-F_n(y)|<\epsilon$.  We can define
equicontinuity at $\infty$ by replacing the condition involving $\delta$
with $y>1/\delta$.  Similarly at $-\infty$.  The point is that one
$\delta$ works for all $n\in\N$.  If $\{F_n\}$ is equicontinuous
at each point of $\Rbar$ we say this sequence is equicontinuous on $\Rbar$.
\begin{corollary}[\cite{ang}, Corollary~3]\label{angequicontinuous}
Let $\{f_n\}\subset\Alex$ such that $f_n\to f$ weakly in $\D$ for some
$f\in\Dp$.  
Suppose $\{F_n\}$ is equicontinuous on $\Rbar$.
Then $f\in\Alex$ and $\|f_n-f\|\to 0$.
\end{corollary}

The proof depends on the Arzel\`a--Ascoli theorem.  See \cite{ang}.

\begin{example}\label{examplesequence1}
{\rm
Let $\{a_n\}$ be a sequence of positive
real numbers that increases to infinity.
Define $f_n$ as the step function
$$
f_n(t)=\left\{\begin{array}{cl}
0, & t\leq n-1 \\
a_n, & n-1<t\leq n\\
-a_n, & n<t\leq n+1\\
0, & t>n+1.
\end{array}
\right.
$$
Then $f_n\in\Alex$ for each $n\in\N$  and  $F_n$ is the piecewise linear
function
$$
F_n(x)=\left\{\begin{array}{cl}
0, & x\leq n-1 \\
a_n(x-n+1), & n-1\leq x\leq n\\
a_n(n+1-x), & n\leq x\leq n+1\\
0, & x\geq n+1.
\end{array}
\right.
$$
It follows that $\|F_n\|_\infty=a_n$.  Note that $F_n\to0$ on $\Rbar$
and that the convergence is quasi-uniform but not uniform.  To see that
it is not uniform, notice that $F_n(n)=a_n\to\infty$.  By 
Theorem~\ref{quasiuniform}, $\intinf f_n\to 0$.
Note that $\{F_n\}$ is uniformly bounded on compact
intervals: $\|F_n\chi_{[a,b]}\|_\infty\leq a_m$ where $m$ is the
largest integer such that $a-1\leq m\leq b+1$.  Hence, $f_n$ converges
weakly to 0 in $\D$.  Theorem~\ref{ang8} and Corollary~\ref{ang9} allow
us to conclude that $\intinf f_n\to 0$.  Also, $\{F_n\}$ is equicontinuous
on $\R$ but not at $\infty$, since if $\delta>0$ then  for integer
$n>1/\delta$ we have $F_n(n)=a_n$ and this can be made arbitrarily
large by taking $n$ large enough.  Hence, Corollary~\ref{angequicontinuous} 
is not applicable.

Although $f_n\to 0$ weakly in $\D$, $\{f_n\}$ does not converge
weakly in $\bv$.  Define $g=\sum_nb_n\chi_{[2n-1,2n]}$ where $\{b_n\}$ is
a sequence of positive real numbers.  Then $Vg =2\sum_n b_n$.
We have $\langle f_{2n},g\rangle=\int_{2n-1}^{2n+1} f_{2n} g =a_{2n}b_n$.  If
$a_n=n^3$ and $b_n=1/n^2$ then $g\in\bv$ but $\langle f_{2n},g\rangle=8n\to
\infty$.

Each function $f_n$ is Riemann integrable and $f_n\to 0$ pointwise on
$\R$ but the sequence of integrals $\intinf f_n$ does not converge
uniformly so the usual convergence theorems for Riemann integration
do not apply.

Convergence theorems for Lebesgue integration also do not apply,
even though each function $f_n\in L^1$.
There is no $L^1$ function that dominates $|f_n|$ for all $n\in\N$
so the dominated convergence theorem is not applicable.  The Vitali
convergence theorem \cite{dunford} gives necessary and sufficient
conditions for taking limits under Lebesgue integrals but is also
not applicable here since $\intinf|f_n|=2a_n\to\infty$, even
though $\intinf f_n=0$ for each $n\in\N$.\qed
}
\end{example}

\begin{example}
{\rm
Let $\{a_n\}$ be a sequence of positive
real numbers such that $a_n/n$ increases to infinity.
Define $f_n$ as the step function
$$
f_n(t)=\left\{\begin{array}{cl}
0, & t\leq 0 \\
a_n, & 0<t\leq 1/n\\
-a_n, & 1/n<t\leq 2/n\\
0, & t>2/n.
\end{array}
\right.
$$
Then $f_n\in\Alex$ for each $n\in\N$  and  $F_n$ is the piecewise linear
function
$$
F_n(x)=\left\{\begin{array}{cl}
0, & x\leq 0 \\
a_nx, & 0\leq x\leq 1/n\\
a_n(\frac{2}{n}-x), & 1/n\leq x\leq 2/n\\
0, & x\geq 2/n.
\end{array}
\right.
$$
It follows that $\|F_n\|_\infty=a_n/n$.  Note that $F_n\to0$ on $\Rbar$
and that the convergence is quasi-uniform but not uniform, since
$F_n(1/n)=a_n/n\to\infty$.  By
Theorem~\ref{quasiuniform}, $\intinf f_n\to 0$.
Note that $\{F_n\}$ is not uniformly bounded on $[0,1]$.
Theorem~\ref{ang8} and Corollary~\ref{ang9} are not applicable.
Also, $F_n$ is not equicontinuous on $[0,1]$ so 
Corollary~\ref{angequicontinuous} is not applicable.  As with 
Example~\ref{examplesequence1}, convergence theorems for Riemann and Lebesgue
integration are not useful here. \qed
}
\end{example}

With Lebesgue integration, the dominated convergence theorem is
particularly useful because it is often easy to find an integrable
function
that dominates each function in a sequence of functions.  There
is a notion of ordering in $\Alex$ that permits monotone and
dominated convergence theorems.  If $f$ and $g$ are in $\Alex$
then $f\geq g$ if $\langle f,\phi\rangle\geq \langle g,\phi\rangle$
for all $\phi\in\D$ such that $\phi\geq 0$.  Then $f\geq g$ if 
and only if $f-g\geq 0$.  It is known that if $f\in\Dp$ and
$f\geq 0$ then $f$ is a Radon measure, i.e., a Borel
measure that is inner and outer regular, and is finite on compact
sets.  
See \cite{ang} for convergence theorems based on this ordering.
A different ordering, more compatible with the Alexiewicz norm, is
described in Section~\ref{sectionlattice} below.

Instead of dominated convergence we have the following convergence
theorem.  We will see in the next section that it is quite useful.

\begin{theorem}\label{theoremconvergencebv}
Let $f\in\Alex$.  Suppose $\{g_n\}\subset\bv$ such that there is
$M\in\R$ so that for all
$n\in\N$, $Vg_n\leq M$.  If $g_n\to g$ on $\Rbar$ for a function $g\in\bv$
then $\lim_{n\to\infty}\intinf fg_n = \intinf fg$.
\end{theorem}
The theorem  is based on Helly's theorem for Riemann-Stieltjes integrals.
See \cite{talvilarae} for a proof.  This paper also contains convergence
theorems for products $f_n g_n$ when $f_n$ is Henstock--Kurzweil integrable.
The proofs carry over to $\Alex$ with no change. 

\section{The Poisson integral and Laplace transform}

A common use of integrals is the integration of functions from a certain
class against a fixed kernel.   We will look at two typical cases, the
Poisson integral 
and Laplace transform.

The upper half plane Poisson integral is given by the convolution
$
u(x,y)=K(x-\cdot,y)\ast f=\intinf f(t)K(x-t,y)\,dt
$,
where the Poisson kernel is $K(x,y)=y/[\pi(x^2+y^2)]$.
It is known that if $f\in L^p$ ($1\leq p\leq \infty$) then $u$ is
harmonic in the upper half plane.  This is also true in $\Alex$.
Fix $x\in\R$ and $y>0$.
Let $f\in\Alex$.  The kernel $t\mapsto K(x-t,y)$ is of bounded variation
on $\Rbar$.  Therefore, the product $f(\cdot)K(x-\cdot,y)$ is in $\Alex$
and $u$ exists on the upper half plane.  To show that we can differentiate
under the integral sign, let $h$ be a nonzero real number and consider 
$$
t\mapsto \frac{K(x+h-t,y)-K(x-t,y)}{h} 
=\frac{-y(2x -2t +h)}{\pi\left[(x+h-t)^2+y^2\right][(x-t)^2+y^2]}.
$$  This function is of bounded variation on $\Rbar$, uniformly for
$h\not=0$.  Hence, using Theorem~\ref{theoremconvergencebv}, we
can differentiate under the integral sign to get
$u_1(x,y)=-\frac{2y}{\pi}\intinf
\frac{f(t)(x-t)\,dt}{\left[(x-t)^2+y^2\right]^2}$.
Similarly, $u_2(x,y)=\frac{1}{\pi}\intinf\frac{f(t)[(x-t)^2-y^2]dt}{
[(x-t)^2 +y^2]^2}$.  And, using these two new kernels and Theorem~\ref{theoremconvergencebv}, we see that $\Delta u(x,y)=\intinf f(t)\Delta K(x-t,y)\,dt=0$
and $u$ is harmonic in the upper half plane.

Using our change of variables Theorem~\ref{changeofvariablestheorem}
with $G(t)=x-t$, $a=-\infty$ and $b=\infty$,  we
can show that $u(x,y)=f(x-\cdot)\ast K(\cdot,y)=\intinf f(x-t)K(t,y)\,dt$.
It is also possible to show that boundary conditions are taken on in
the Alexiewicz norm, i.e., $\|u(\cdot,y)-f(\cdot)\|\to 0$ as $y\to 0^+$.

Let $\R^+=(0,\infty)$
and $f\in\Alex(\R^+)$.  We will say that the variation of a complex-valued
function is the sum of the variations of the real and imaginary parts.
Let $x,y\in\R$ and write
$z=x+iy$.
The function $t\mapsto e^{-zt}$ is of bounded
variation on $[0,\infty]$ if $x>0$ or if $z=0$. Hence, the Laplace transform
of $f$ is
$\hat{f}(z)=\int_0^\infty f(t)\,e^{-zt}\,dt$ and exists for $x>0$ or $z=0$.   We can now prove some basic
properties of the Laplace transform.
First we will prove $\fhat$ is 
differentiable.  Fix $x>0$ and take $h\in\C$ such that $0<|h|< x/2$.
For fixed $z=x+iy$ with $x>0$ write
$g_h(t)=[\exp(-(z+h)t)-\exp(-zt)]/h$.
Then
$$
|g_h'(t)|  =  e^{-xt}\left|\frac{(z+h)e^{-ht}-z}{h}\right| =
e^{-xt}\left|\left(e^{-ht}-1\right)\frac{z}{h} + e^{-ht}\right|.
$$
By Cauchy's theorem, 
$$
e^{-ht}=1+\frac{h}{2\pi i}\int_C\frac{e^{-st}\,ds}{s(s-h)}
$$
where $C$ is the circle with centre $0$ and radius $x/2$ in the complex plane.  Then
$|g_h'(t)|\leq (2|z|/[x-2|h|]+1)e^{-xt/2}$ and 
$Vg_h\leq (2|z|/[x-2|h|]+1)(2/x)$ so that $g_h$ is of bounded variation
on $[0,\infty]$, uniformly as $h\to 0$.  By Theorem~\ref{theoremconvergencebv},
$d\fhat(z)/dz=
-\int_0^\infty f(t)\,te^{-zt}\,dt$.   Similarly, we can differentiate
under the integral sign to get $d^n\fhat(z)/dz^n=(-1)^n\int_0^\infty
f(t)\, t^n e^{-zt}\,dt$ for all $n\in\N$.

One difference between Laplace transforms in $\Alex(\R^+)$ 
and Laplace transforms
of distributions is that we get
a different growth condition as $z\to\infty$.
Write $z=x+iy$ with $x>0, y\in\R$.
Let $\delta>0$.
Integrate by parts to get $\fhat(z)=z\int_0^\delta F(t)e^{-zt}\,dt
+z\int_\delta^\infty F(t)e^{-zt}\,dt$.  Then
$|\fhat(z)|\leq (|z|/x)\max_{[0,\delta]}|F|+(|z|/x)\|F\|_\infty e^{-x\delta}$.  Given
$\epsilon>0$, take $\delta$ small enough so that 
$\max_{[0,\delta]}|F|<\epsilon$. Let $0\leq\alpha<\pi/2$.  
We then have $\fhat(z)=o(1)$ as $z\to\infty$ in the cone $|\arg(z)|\leq\alpha$.
We can show this estimate is sharp by showing it is sharp as $z=x$ goes to
infinity on the positive real axis.
Suppose $A\fn(0,\infty)\to(0,1)$ with $\lim_\infty
A=0$.  First show $A$ has a suitably smooth majorant.  Define
$B(s)=\sup_{0<t\leq s}eA(1/t)$.  Then $B(s)\geq eA(1/s)$ for all $s>0$,
$B$ is increasing and $\lim_{s\to 0^+}B(s)=0$.  Now define
$$
F(s)  =  \left\{\!\!\!\!\!\!\begin{array}{cl}
\left[B\left(\frac{1}{n}\right)-B\left(\frac{1}{n+1}\right)\right]
(n+1)(n+2)s & \\
\quad-(n+1)B\left(\frac{1}{n}\right) +(n+2)B\left(\frac{1}{n+1}\right), & 
\frac{1}{n+2}\leq s\leq \frac
{1}{n+1}\mbox{ for some } n\in\N\\
0, & s=0\\
B(1), & s\geq 1/2.
\end{array}
\right.
$$
Then $F\in\Ct(\R^+)$, $F(s)\geq eA(1/s)$ for all $s\in(0,1/2]$.  Since $F$ is
increasing and piecewise linear, $F\in AC(\R^+)\cap C^0([0,\infty])$.  
Let $f=F'$ and let $s\in(0,1/2]$.
Then $\fhat(x)\geq \int_0^s f(t)e^{-xt}\,dt\geq F(s)e^{-xs}$.  Now suppose
$x \geq 2$.  Let $s=1/x$.  Then $\fhat(x)\geq F(1/x)e^{-1}\geq A(x)$.  Hence,
the estimate $\fhat(z)=o(1)$ ($z\to\infty$, $|\arg(z)|\leq\alpha$) is sharp, not only in $\Alex$
but in $L^1$ as well.  Note that for the Dirac distribution,
$\widehat{\delta}(z)=\exp(0)=1$ so the
estimate does not hold for measures or distributions that are the second
derivative of a continuous function.  For distributions in general, the
Laplace transform can have polynomial growth.  See \cite[p.~236, 237]{zemanian}.

Since the kernel decays
exponentially, we can define a Laplace transform under
weaker conditions.  Define the locally integrable distributions on $[0,\infty)$
by $\Alex(loc)=\{f\in\Dp(\R^+)\mid
 f=F' \mbox{ for some } F\in C^0([0,\infty))\}$.  In this case, $f=F'$ means
that for all $\phi\in\D(\R^+)$ we have $\langle f,\phi\rangle=
-\langle F,\phi'\rangle$.  For $f\in\Alex(loc)$ there is a continuous
function $F$ such that $\int_0^x f= F(x)-F(0)$ for all $x\in[0,\infty)$.
Note that $\lim_\infty F$ need not exist.
Let $r\in\R$.  Define $F_r(x)=\int_0^x f(t)e^{-rt}\,dt =F(x)e^{-rx}-F(0)
+r\int_0^x F(t) e^{-rt}\,dt$.  Note that $F_r(0)=0$ and 
$F_r\in C^0([0,\infty))$.  Now we can define the weighted space
$\Alex[e^{r\cdot}]=\{f\in\Alex(loc)\mid f=F' \mbox{ for some } 
F\in C^0([0,\infty)) \mbox{ such that } \lim_\infty F_r \mbox{ exists in }
\R\}$.  For example, if $F$ is a continuous function such that
$F(x)e^{rx}/x^2$ is bounded as $x\to\infty$ then $F\in\Alex[e^{r\cdot}]$.
We then have $\int_0^\infty f(t)e^{-rt}\,dt=\lim_{x\to\infty}F_r(x)$.
The limit is independent of which primitive $F\in C^0([0,\infty))$ is
used. If $f\in\Alex[e^{r\cdot}]$ then ${\hat f}(z)$ exists for all
$z\in\C$ such that
${\cal R}e(z)> r$ or ${\cal R}e(z)\geq r, {\cal I}m(z)=0$.  If $f$ is in one of these exponentially weighted spaces there
are similar differentiation and growth results as to when $f\in\Alex(\R^+)$.
Using an analogous technique, we can define weighted integrals $\intinf fg$
for functions $g$ 
that are of locally bounded variation.

\section{Banach lattice}\label{sectionlattice}

In $\Ct$ there is the pointwise order: for $F,G\in\Ct$, $F\leq G$ if and
only if $F(x)\leq G(x)$ for all $x\in\Rbar$.  It is easy to see that this
relation is {\it reflexive} ($F\leq F$), {\it antisymmetric}
($F\leq G$ and $G\leq F$ imply $F=G$), and {\it transitive} ($F\leq G$
and $G\leq H$ imply $F\leq H$).  This puts a {\it partial order} on $\Ct$.

As $\Alex$ is isomorphic to $\Ct$, it inherits this partial order.
For $f,g\in\Alex$, we define $f\leq g$ if and only if $F\leq G$.
For example,
let $f(t)=\sin(t)/t$ for $t>0$ and $f(t)=0$ for $t<0$.  Then $f\in\Alex$.
We have $F(x)=\int_0^x f$ for $x\geq 0$ and $F(x)=0$ for $x\leq 0$.  This
is the {\it sine integral}, ${\rm Si}(x)$, and it is easy to show $F(x)\geq 0$
for all $x\in\R$.  Hence, $f\geq 0$ in $\Alex$.  This ordering on $\Alex$ is then not
compatible with the usual pointwise ordering that we can use in $L^1$, i.e.,
$f\geq g$ if and only if $f(t)\geq g(t)$ for almost all $t\in\R$. The
function defined by
$\max(f(t),0)$ is not in $\Alex$.  Nor is our ordering compatible with
the usual one for distributions: if $T\in\Dp $ then $T\geq 0$ if 
and only if $T$ is a 
Radon measure.  The function $f(t)=\sin(t)/t$ is not positive in the 
distributional sense.  It is not even the difference of two positive,
Lebesgue integrable functions so it is not a signed measure.  In $\Alex$,
the relation $f\geq 0$ means that for each $x\in\R$, the integral
over $(-\infty, x]$ is not negative, i.e., to the left of $x$ there is more
positive stuff than negative stuff.  It is a not a {\it linear ordering}.
For example, $f(t)=-2t\exp(-t^2)$ and $g(t)=-2(t-1)\exp(-(t-1)^2)$ are not comparable.

Now, $\Ct$ is closed under the operations
$(F\vee G)(x)=\sup(F(x),G(x))=\max(F(x),G(x))$ and
$(F\wedge G)(x)=\inf(F(x),G(x))=\min(F(x),G(x))$.
It is then a {\it lattice}.
And, $\Ct$ is also a {\it Banach lattice}.  This means that the order is
compatible with the vector space operations and norm.  For all $F,G\in\Ct$,
\begin{enumerate}
\item[(i)]
$F\leq G$ implies
$F+H\leq G+H$ for all $H\in\Ct$
\item[(ii)]
if $F\leq G$ then $aF\leq aG$ for all
real numbers $a\geq 0$
\item[(iii)] $|F|\leq |G|$ implies $\|F\|_\infty\leq \|G\|_\infty$.
\end{enumerate}
A good introduction to lattices can be found in \cite{aliprantis}.

As usual, in $\Ct$ we define $F^+=F\vee 0$, $F^-=F\wedge 0$ and
$|F|=F\vee (-F)$.  The
{\it Jordan decomposition} is $F=F^+-F^-$.  It is also true that
$|F|=F^+ +F^-$.  In $\Alex$, $f^+=(F^+)'$, $f^-=(F^-)'$ and $|f|=|F|'$.
These definitions make sense since $F\in\Ct$ so $F^+$, $F^-$ and
$|F|$ are all in $\Ct$ and then their derivatives are in $\Alex$.  For the
function $f(t)=\sin(t)/t$ when $t>0$ and $f(t)=0$, otherwise, we have
$f^+=|f|=f$ and $f^-=0$.

\begin{theorem}\label{theoremlattice}
$\Alex$ is a Banach lattice.
\end{theorem}

\noindent
{\bf Proof:} First we need to show that
$\Alex$ is closed under the operations $f\vee g$ and
$f\wedge g$.  For $f,g\in\Alex$, we have
$f\vee g  =  \sup(f,g)$. This is $h$ such that
$h\geq f$,  $h \geq g$, and if $h_1\geq f$, $h_1\geq g$, then $h_1\geq h$.
This last statement is equivalent to 
$H\geq F$,  $H \geq G$, and if $H_1\geq F$, $H_1\geq G$, then $H_1\geq H$.
But then $H=\max(F,G)$ and $h=H'$ so $f\vee g=(F\vee G)'\in\Alex$.
Similarly, $f\wedge g=(F\wedge G)'\in\Alex$.

If $f,g\in\Alex$ and $f\leq g$ then $F\leq G$.  Let $h\in\Alex$.  Then, $F+H\leq G+H$.
But then $(F+H)'=F'+H'=f+h\leq g+h$.  If  $a\in\R$ and $a\geq 0$
then $(aF)'=aF'=af$ so $af\leq ag$.  And, if $|f|\leq |g|$ then 
$|F|'\leq |G|'$ so $|F|\leq |G|$, i.e., $F(x)\leq G(x)$ for all $x\in\Rbar$.
Then $\|f\|=\|F\|_\infty\leq \|G\|_\infty=\|g\|$. And, $\Alex$ is a
Banach lattice that is isomorphic to $\Ct$.\qed

Linearity of the derivative was necessary to prove conditions (i) and (ii),
whereas, for (iii) we needed the
fact that $\Ct$ and $\Alex$ are isometric.  It is a fact that every 
Banach lattice
is isomorphic to the vector space of continuous functions on some compact
Hausdorff space.  See, for example, \cite[pp.~395]{dunford}.

The following results follow immediately from the definitions.

\begin{theorem}
Let $f,g\in\Alex$.
{\rm (a)} If $f\leq g$ then $F(x)\leq G(x)$ for all $x\in\Rbar$.
{\rm (b)} If $\int_{-\infty}^x f\leq \int_{-\infty}^x g$ for all
$x\in\R$ then $f\leq g$.
{\rm (c)} $|f|\in\Alex$ and $|\int_{-\infty}^x f|\leq \int_{-\infty}^x|f|$ for all $x\in\Rbar$.
{\rm (d)} $\|\,|f|\,\|=\|\,|F'|\,\|=\|\,|F|\,\|_\infty= \|f\|$.
\end{theorem}

The order on $\Alex$ gives us absolute integration since if $F$ is
continuous, so is $|F|$ and then integrability of $f$ implies integrability
of $|f|$.
Notice that the definition of order allows us to integrate both sides of
$f\leq g$ in
$\Alex$ to get $F\leq G$ in $\Ct$.  The isomorphism allows us to differentiate
both sides of 
$F\leq G$ in $\Ct$ to get $F'\leq G'$ in $\Alex$.  However, there is no
pointwise implication.  For example, 
$F(x)\geq 0$ for all $x\in\R$ does not imply $F'(x)\geq 0$ for all $x\in\R$. 
Take
$F(x)=\exp(-x^2)$.  
And, if $f$ and $g$ are functions
in $\Alex$ and $f(t)\leq g(t)$ for all $t\in\R$, we cannot conclude
that $f\leq g$ in 
$\Alex$.  This was shown with the $f(t)=\sin(t)/t$ function above.
Note also that the partial ordering mentioned at the end of 
Section~\ref{sectionconvergence} fails to be a vector lattice.  If $f\in\Alex$
is a function and $\langle f,\phi\rangle \geq 0$ for all $\phi\in\D$ with
$\phi\geq 0$ then $f\geq 0$ almost everywhere.  Hence, $\sup(f,0)$ need not
be in $\Alex$.  This is the case for any function that has a conditionally
convergent integral, as with our $\sin(t)/t$ function.  
In the next section we consider the more usual type of absolute integrability.

\section{Absolute convergence}

Suppose $f\in\Alex$.
Let $\|f\|_{{\cal ABS}}=\sup_{\stackrel{\phi\in\D}{\|\phi\|_\infty\leq 1}}
\langle f,\phi\rangle$ and define
${\cal ABS}=\{f\in\Alex\mid\|f\|_{{\cal ABS}}<\infty\}$.  We will show
that ${\cal ABS}$ provides a sensible extension of the notion of absolute
integrability.  If
$f\in\Alex$ and its primitive is $F\in\bv\cap\Ct$ then, by the
H\"older inequality,
$$
|\langle f,\phi\rangle| = |\langle F',\phi\rangle|=\left|\intinf F'\phi\right|
\leq 2VF\,\|\phi\|_\infty.
$$
So, $f\in{\cal ABS}$.  If $f\in{\cal ABS}$ then
$$
\sup_{\stackrel{\phi\in\D}{\|\phi\|_\infty\leq 1}}\langle f,\phi\rangle
=\sup_{\stackrel{\phi\in\D}{\|\phi\|_\infty\leq 1}}\intinf F\phi'<\infty.
$$
Since $F\in\Ct$ we have
$VF=essvar \, F<\infty$.
Thus, $f\in{\cal ABS}$ if and only if $VF<\infty$.  
See Section~\ref{integrationbyparts} for the definition of the essential
variation.

From the
definition of variation it follows that $\|f\|_{{\cal ABS}}=VF$.
We know $\bv$ is a Banach space.  Clearly $\Ct\cap \bv$ is a subspace.
To show it is complete, suppose $\{F_n\}\subset \Ct\cap \bv$ is Cauchy
in the $\bv$ norm.  Then there is $F\in\bv$ such that $V(F_n-F)\to 0$.
We need to show $F\in\Ct$. Let $x\in\Rbar$.  We have 
\begin{eqnarray*}
|F(x)-F(y)| & \leq & |F(x)-F_n(x)-F(y)+F_n(y)|+ |F_n(x)-F_n(y)|\\
 & \leq & V(F_n-F)+|F_n(x)-F_n(y)|.
\end{eqnarray*}
Given $\epsilon>0$ we can take $n$ large enough so that $V(F_n-F)<\epsilon/2$.
Since $F_n\in\Ct$ we can now take $y$ close enough to $x$ so that
$|F_n(x)-F_n(y)|<\epsilon/2$.  Hence, $F\in\Ct$ and $\Ct\cap \bv$ is
a Banach space.  The integral provides a linear isometry between
${\cal ABS}$ and $\Ct\cap \bv$.  Hence, $\|f\|_{{\cal ABS}}$ is a norm and
${\cal ABS}$ is a Banach space.  We identify ${\cal ABS}$ as the subspace
of $\Alex$ consisting of absolutely integrable distributions by analogue
with the fact that primitives of Denjoy or wide Denjoy integrable functions
need not be of bounded variation but primitives of $L^1$ functions are
absolutely continuous and hence of bounded variation.

\section{Odds and ends}\label{odds}

We collect here various other results.  The first is that there are no
improper integrals. 

\begin{theorem}[Hake Theorem]\label{haketheorem}
Suppose $f\in\Dp$ and $f=F'$ for some $F\in C^0(\R)$.  If $\lim _{\infty}F$
and $\lim_{-\infty}F$ exist in $\R$ then $f\in\Alex$ and
$\intinf f=\lim_{x\to\infty}\int_{0}^x f +\lim_{x\to-\infty}\int_{x}^0 f$.
\end{theorem}

\bigskip
\noindent
{\bf Proof:} Define ${\overline F}(x)=F(x)$ for $x\in\R$, 
${\overline F}(\infty)=\lim_\infty F$, 
${\overline F}(-\infty)=\lim_{-\infty} F$.  Then ${\overline F}\in C^0(\Rbar)$
and ${\overline F}\,'=f$.  Hence, $f\in\Alex$ and
\begin{eqnarray*}
\intinf f & = & {\overline F}(\infty) - {\overline F}(-\infty)\\
 & = & \lim_{\infty}F-\lim_{-\infty}F\\
 & = & \lim_{x\to\infty}\left[F(x)-F(0)\right]
+\lim_{x\to-\infty}\left[F(0)-F(x)\right].\qed
\end{eqnarray*}

There are similar versions on compact intervals and intervals such as
$[0,\infty)$.  The corresponding result is false for Lebesgue integrals.
For example, $\lim_{x\to\infty}\int_0^x \sin(t^2)\,dt={\sqrt \pi}/(2^{3/2})$,
but the function $t\mapsto \sin(t^2)$ is not in $L^1$.  The integral
is called a {\it Cauchy--Lebesgue} integral and in this case is also
an improper Riemann integral.  The theorem is
true for Henstock--Kurzweil integrals.
Proving the Hake
theorem
for the Henstock--Kurzweil or Perron integral is more involved.
See \cite{gordon},
Theorem~9.21 and
Theorem~8.18.

\begin{theorem}[Second mean value theorem]\label{secondmeanvaluetheorem}
Let $f\in\Alex$ and let $g\fn\Rbar\to\R$ be monotonic. Then
$\intinf fg= g(-\infty)\int_{-\infty}^\xi f+ g(\infty)\int_{\xi}^\infty f$
for some $\xi\in\Rbar$.
\end{theorem}

\bigskip
\noindent
{\bf Proof:} Integrate by parts and use the mean value theorem for
Riemann--Stieltjes integrals \cite[\S7.10]{mcleod}:
\begin{eqnarray*}
\intinf fg & = & F(\infty)g(\infty)-\intinf F\,dg\\
 & = & F(\infty)g(\infty)-F(\xi)\intinf dg\\
 & = &  F(\infty)g(\infty)-F(\xi)[g(\infty)-g(-\infty)]\\
 & = & g(-\infty)F(\xi) + g(\infty)[F(\infty)-F(\xi)].\qed
\end{eqnarray*}
This proof is taken from \cite{celidze}, where a proof of the Bonnet
form of the second mean value theorem can also be found.

Using the distributional integral, it is possible to formulate
a version of Taylor's theorem with integral remainder.  For
an approximation by an $n$th degree polynomial it is only required
that $f^{(n)}$ be continuous.

\begin{theorem}[Taylor]
Suppose $[a,b]\subset \R$.
Let $f\fn[a,b]\to\R$ and let $n\geq 0$ be an integer.  If $f^{(n)}\in
C^0([a,b])$ then for all $x\in[a,b]$ we have $f(x) =P_n(x) + R_n(x)$ where
$$
P_n(x)=\lsum_{k=0}^n\frac{f^{(k)}(a)(x-a)^k}{k!}
$$
and
$$
R_n(x)=\frac{1}{n!}\int_a^x f^{(n+1)}(t)(x-t)^{n}\,dt.
$$
For each $x\in [a, b]$ we have the estimate
\begin{eqnarray*}
|R_n(x)| & \leq & \frac{(x-a)^n\|f^{(n+1)}\chi_{[a,x]}\|}{n!}\leq
\frac{(x-a)^n\|f^{(n+1)}\|}{n!}\\
  & = & 
\frac{(x-a)^n}{n!}\max_{a\leq\xi\leq x}|f^{(n)}(\xi)-f^{(n)}(a)|.
\end{eqnarray*}
And,
$$
\|R_n\|\leq \|R_n\|_1\leq \frac{(b-a)^{n+1}}{(n+1)!}\|f^{(n+1)}\|=\frac{(b-a)^{n+1}}{(n+1)!}
\max_{a\leq \xi\leq b}|f^{(n)}(\xi)-f^{(n)}(a)|.
$$
\end{theorem}

The remainder exists since
the function $t\mapsto (x-t)^n$ is monotonic for each $x$. 
Repeated integration by parts establishes the integral remainder
formula.
Estimates of the remainder follow upon applying the second mean value
theorem.
See \cite{talvilataylor} for various other estimates of the
remainder. 
Usual versions
of Taylor's theorem require $f^{(n+1)}$ to be integrable.  For the Lebesgue
integral this means taking $f^{(n)}$ to be absolutely continuous.  
Here we only need
$f^{(n)}$ continuous.

\begin{theorem}[Homogeneity of Alexiewicz norm]\label{theoremhomogenous}
Let $f\in\Alex$. 
For $t\in\R$,
define the translation $\tau_t$ by 
$\langle\tau_tf, \phi\rangle=\langle f, \tau_{-t}\phi\rangle$
where $\tau_t\phi(x)=\phi(x-t)$ for
$\phi\in\D$.
The Alexiewicz norm is translation invariant: If $f\in\Alex$ then
$\tau_tf \in\Alex$
and $\|\tau_tf\|=\|f\|$.  Translation is continuous:  $\|f-\tau_tf\|\to
0$ as $t\to 0$.
\end{theorem}

\bigskip
\noindent
{\bf Proof:}
If $f\in\Alex$ then
a change of variables shows 
\begin{eqnarray*}
\langle \tau_tf,\phi\rangle & = & \langle f,\tau_{-t}\phi\rangle
  = 
\intinf f(s)\phi(s+t)\,ds
  =  \intinf f(s-t)\phi(s)\,ds\\
  & = &  
\intinf F'\!(s-t)\phi(s)\,ds
  =  \intinf(\tau_tF)'\!(s)\phi(s)\,ds
\end{eqnarray*}
 and $\tau_tF\in\Ct$ is the primitive of $\tau_t f$.
Hence, $\tau_tf\in\Alex$.  It is clear that 
$\|F\|_\infty=\|\tau_tF\|_\infty$ for all $t\in\R$.  Hence, $\|f\|=\|\tau_tf\|$.

As well,
\begin{eqnarray*}
\sup_{x\in\R}\left|\int_{-\infty}^x\left[f(s)-\tau_tf(s)\right]ds\right|
 & = & \sup_{x\in\R}\left|F(x)-F(x-t)\right|\\
 & \to & 0 \mbox{ as } t\to 0\mbox{ since } F \mbox{ is uniformly 
continuous.}\qed
\end{eqnarray*}
See \cite{talvilaalex} for some other continuity properties of the
Alexiewicz norm.

A Banach space satisfying
the conditions of Theorem~\ref{theoremhomogenous} is called
{\it homogeneous}. 

\begin{theorem}[Equivalent norms]\label{Equivalent norms}
The following norms on $\Alex$ are equivalent to $\|\cdot\|$.  For
$f\in\Alex$, define $\|f\|'=
\sup_I|\int_If|$ where the supremum is taken over all compact intervals 
$I\subset\R$;
$\|f\|''=\sup_{g}\int fg$, where the supremum is taken over all
$g\in\bv$ such that $|g|\leq 1$ and $Vg\leq 1$; $\|f\|'''=\sup_{g}\int fg$, 
where the supremum is taken over all
$g\in\ebv$ such that $\|g\|_\infty\leq 1$ and $essvar \, g\leq 1$.
\end{theorem}

\bigskip
\noindent
{\bf Proof:} We have $\|f\|'=\sup_{a<b}\left|\int_a^b f\right|=
\sup_{a<b}\left| F(b)-F(a)\right|\leq 2\|f\|$.  And, $\|f\|\leq \|f\|'$.
Hence, $\|\cdot\|$ and $\|\cdot\|'$ are equivalent.  Let $g\in\bv$
with $|g|\leq 1$ and $Vg\leq 1$. 
By the H\"older
inequality (Theorem~\ref{holder}),
$$
\left|\intinf fg\right| \leq \|f\|\left[\inf|g|+2Vg\right]\leq 3\|f\|.
$$
And,
$$
\|f\|''\geq \max\left(\sup_{x\in\R}\intinf f\chi_{(-\infty,x]},
-\sup_{x\in\R}\intinf f\chi_{(-\infty,x]}\right).
$$
It follows that $\frac{1}{3}\|f\|''\leq \|f\|\leq \|f\|''$.  
The proof for $\|\cdot\|'''$ is similar.\qed

The following definition allows us to integrate any distribution over a compact interval.  The result is also a distribution.  If $T\in\Dp$ and 
$[a,b]\subset\R$, define 
\begin{eqnarray*}
\left\langle\int_a^b T',\phi\right\rangle & := &
\left\langle T', \int_a^b\tau_t\phi(\cdot)\,dt\right\rangle
= -\langle T, \int_a^b\phi'(\cdot -t)\,dt\rangle\\
 & = &  \langle T, \tau_b\phi\rangle-\langle T,\tau_a\phi\rangle =
 \langle\tau_{-b}T,\phi\rangle-
\langle\tau_{-a}T,\phi\rangle.
\end{eqnarray*}
The translation $\tau_{-a}$ was defined in Theorem~\ref{theoremhomogenous}.
In the case of $T'=f\in\Alex$ this gives 
$
\left\langle \int_a^b f,\phi\right\rangle  =  \intinf F(t) \phi(t-b)\,dt - 
\intinf F(t) \phi(t-a)\,dt$,
which is a convolution.  Since $F$ is continuous, we can recover the
value $\int_a^b f\in\R$ by evaluating on a delta sequence $\{\phi_n\}$.
See the end of Section~\ref{distributions}.  We then have
$\langle\int_a^b f,\phi_n\rangle\to F(b)-F(a)$.  This method of integration was
developed by J.~Mikusi\'nski, J.A.~Musielak and R.~Sikorski
 in the 1950's and 1960's
\cite{jmikusinski}, \cite{musielak}, \cite{sikorski}.  The advantage is that
it can integrate every distribution over a compact interval.  The disadvantage 
is that integrals over $(-\infty, \infty)$ must be treated as improper 
integrals since $\tau_{\pm\infty}\phi=0$.  As we saw in 
Theorem~\ref{haketheorem}, there are no improper integrals 
in $\Alex$.  And, of course $\Alex$ is a Banach space, whereas $\Dp$ is not.
 
\section{Further threads}

In this final section we list several topics in passing and several
ideas for further research.

\noindent
{\bf 1. What happened to the measure?}  In Lebesgue and Henstock--Kurzweil
integration the measure appears explicitly.  With  the distributional 
integral it is disguised in the
formula $F'=f$, out of which $\langle f, \phi\rangle =-\langle F,\phi'\rangle$
for all $\phi\in\D$.  The derivative is
$$
\phi'(x)  =  \lim_{h\to 0}\frac{\phi(x+h)-\phi(x)}{h}
  =  \lim_{h\to 0^+} \frac{\phi(I(x,h))}{\lambda(I(x,h))},
$$
where $I(x,h)$ is the interval centred on $x$ with radius $h$ and
we have replaced $\phi$ by the interval function $\phi((a,b))=\phi(b)-\phi(a)$.
Replacing Lebesgue measure $\lambda$ with some other measure $\mu$ gives
the Radon--Nikodym derivative with respect to $\mu$.  To integrate $f$ with
respect to $\mu$ we need to use the Radon--Nikodym derivative when we define
integration by parts for distributions.  The test functions would have
to have all their Radon--Nikodym derivatives continuous with respect to $\mu$.
The primitives would have to be continuous with respect to $\mu$, rather 
than pointwise.  For continuity at $x$ this means that  
for all $\epsilon>0$ there is $\delta>0$ such that
$\mu(I(x,|x-y|))<\delta$ gives $|F(x)-F(y)|<\epsilon$, whereas replacing
$\mu$ with $\lambda$ gives the usual pointwise definition of continuity.\\

\noindent
{\bf 2. Integration in $\R^n$.}  The Denjoy integral has not been easy to formulate in $\R^n$ due to 
the difficulty
of defining $ACG*$ in $\R^n$.  For the fearless, see Chapter~2 in 
\cite{celidze}.   There is, however, a distributional integral in $\R^n$.
If $f\in\Dp(\R^n)$ then $f$ is integrable
if there is a function $F\in C^0(\Rbar^n)$ such that $DF=f$.  The
differential operator is $D=\frac{\partial^n}{\partial x_1\partial x_2 \cdots
\partial x_n}$.  Now,
$\langle f,\phi\rangle = \langle DF, \phi\rangle=(-1)^n\langle F,D\phi\rangle$
where $\phi$ is a $C^\infty$ function with compact support in $\R^n$.
For example,
$\int_a^b\int_c^d F_{12}=F(b,d)-F(a,d)-F(b,c)+F(a,c)$ for each continuous
function $F$.
This is the form of the integral given in \cite{pmikusinski1}.
For details see \cite{ang}, where there are applications to the wave equation
and theorems of  Fubini and Green.
This definition extends the Lebesgue and Henstock--Kurzweil integrals.  
But, it is not invariant under rotations since the operator $D$ is
not invariant under rotations.  For example, a rotation of $\pi/4$ for
which $(x,y)\mapsto (\xi,\eta)$ transforms
$D$ into the wave operator $\partial^2/\partial\xi^2-\partial^2/\partial\eta^2$.
Hence, if $f$ is integrable its rotation need not be integrable.

W.~Pfeffer \cite{pfeffer} has defined a 
nonabsolute integral 
that is 
invariant under rotations and other transformations but it is based on
different principles.  In some sense, his integral is designed
to invert the divergence operator.  A possible extension of  Pfeffer's integral 
in the spirit of distributional integrals can be obtained with  the following
definitions.
If $g\in L^1_{loc}(\R^n)$ then $g$ is of {\it local bounded variation}
if $\sup\int_U f\,{\rm div}\phi<\infty$ for each open ball $U\subset\R^n$, 
where the supremum is taken over all
$\phi\in{\cal D}(U)$ with $\|\phi\|_\infty\leq 1$.
A measurable set $E\subset\R^n$ has {\it locally
finite perimeter} if $\chi_E$ is of local bounded variation.  
Sets
with Lipshitz boundary have this property and thus polytopes do as well.
Suppose $\Omega\subset\R^n$ is open and $E\subset\Omega$ 
has locally finite perimeter.
Then
$f\in{\cal D}'(\Omega)$ is integrable over $E$ if there is a continuous function
$F\fn\overline{E}\to\R^n$ such that $f={\rm div} F$ in ${\cal D}'(\Omega)$.  Then
$$
\int_E f= \int_E{\rm div} F=\int_{\partial*E}F\cdot n\,
d{\cal H}^{n-1}
$$
where $\partial*E$ is the measure-theoretic boundary of $E$,
$n$ is the outward normal and ${\cal H}^{n-1}$ is Hausdorff measure.
The final integral exists since $F$ is continuous.  This definition of
the integral is based on the Gauss--Green theorem, whose usual version
requires $F$ to be $C^1$.  See \cite{evans} or \cite{ziemer}.

Note that if $F$ is a continuous function in $\R^2$ and $f=F_{21}$ in
$\Dp(\R^2)$ then
$f={\rm div}(F_2, 0)$.  Since the boundary of a Cartesian interval in
$\R^2$ is a union of four intervals in $\R$, the above integral
can be used twice to obtain the formula 
$\int_a^b\int_c^d F_{12}=F(b,d)-F(a,d)-F(b,c)+F(a,c)$.  
Hence, the Gauss--Green integral includes the integral
of Mikusi\'nski and
Ostaszewski \cite{pmikusinski1};
Ang, Schmidt  and Vy \cite{ang}.\\

\noindent
{\bf 3. The regulated primitive integral.}  A function on the real line is {\it regulated} if it has
a left and right limit at each point.  It is known that the Riemann--Stieltjes
integral $\intinf F\,dg$ exists when one of $F$ and $g$ is regulated and
the other 
is of bounded variation.  We can then replace $\Ct$ with the space of
regulated functions.  Then we can integrate all distributions that are the
distributional
derivative of a regulated function.  If $f=F'$ then 
there are four integrals
$\int_{(a,b)}f=F(b-)-F(a+)$, 
$\int_{[a,b)}f=F(b-)-F(a-)$, 
$\int_{(a,b]}f=F(b+)-F(a+)$, 
$\int_{[a,b]}f=F(b+)-F(a-)$,
which need not be same since the left and right limits
of $F$ are not necessarily equal. 
This will allow us to integrate signed Radon measures
since if $\mu$ is a signed Radon measure then $F(x):=\int_{-\infty}^x d\mu$ is a
function of bounded variation and hence regulated.  For example, the Dirac
distribution is the derivative of the Heaviside step function, 
$H(x)=1$ for $x\geq 0$ and $H(x)=0$, otherwise.  And, $\int_{(0,1)} \delta =
\int_{(0,1)} H'=H(1-)-H(0+)= 1-1=0$. Whereas,
 $\int_{[0,1)} \delta =
H(1-)-H(0-)= 1-0=1$.  The regulated primitive integral will be discussed
in detail elsewhere \cite{talvilaregulated}.

It is not clear if we get a useful integral by replacing $\Ct$ with such
Banach spaces as $L^p$ ($1\leq p\leq\infty$) or $\bv$.\\

In light of the existence of other integrals that invert distributional
derivatives, we propose the name {\it continuous primitive integral} for
the integral described in this paper.

\end{document}